\documentclass[12pt]{article}
\usepackage{a4}
\usepackage[ctr1,eng]{myarticl}
\usepackage{amsmath}
\usepackage{amssymb}
\usepackage[frame,curve,knot]{xy}
\usepackage{epsfig}

\newcommand{\mwedge }{{\scriptstyle \land }}
\newcommand{\mc }{\mathcal }

\newcommand{\comp }{\mathbb{C}}
\newcommand{\compx }{\comp ^\times }
\newcommand{\real }{\mathbb{R}}
\newcommand{\realx }{\real ^\times }
\newcommand{\Oq }[1]{\mathrm{O}_{q^2}(#1)}
\newcommand{\LOq }[1]{\Lambda (\mathrm{O}_{q^2}^{#1})}
\newcommand{\OOq }[1]{\mathcal{O}(\mathrm{O}_{q^2}(#1))}

\newcommand{\Uqso }[1]{U_{q^2}(\mathfrak{so}_{#1})}

\newcommand{\Clg }[2][N]{\mathrm{Cl} _{q^2}(#1,#2)}
\newcommand{\Cl }[1][N]{\Clg[#1]{c^2}}
\newcommand{\Cls }[2][N]{\mathrm{Cl}^{#2}_{q^2}(#1,c^2)}
\newcommand{\Cle }[1][N]{\Cls[#1]{+}}
\newcommand{\Clo }[1][N]{\Cls[#1]{-}}
\newcommand{\Clz }[1][N]{\widehat{\mathrm{Cl}}_{q^2}(#1,c^2)}

\newcommand{\qp }{[2]}
\newcommand{\qm }{\hat{q}}

\newcommand{\ot }{\otimes }

\newcommand{\gcl }{\gamma }
\newcommand{\dind }[2]{{\genfrac{}{}{0pt}1{#1}{#2}}}

\newcommand{\im }{\mathrm{im}\,}

\newcommand{\End }[1]{\mathrm{End}(#1)}
\newcommand{\copr }{\Delta }
\newcommand{\coun }{\varepsilon }
\newcommand{\antip }{S}

\newcommand{\ctr }[2]{\bigl\langle \!\!\langle #1,#2\rangle \!\!\bigr\rangle }
\newcommand{\lact }{\triangleright }
\newcommand{\qnum }[1]{\lbrack \! \lbrack #1 \rbrack \! \rbrack }
\newcommand{\qbinom }[2]{\left\lbrack \!\! \left\lbrack
{\genfrac{}{}{0pt}1{#1}{#2}}\right\rbrack \!\! \right\rbrack }
\newcommand{\bigqbinom }[2]{\left\lbrack \!\! \left\lbrack
{\genfrac{}{}{0pt}0{#1}{#2}}\right\rbrack \!\! \right\rbrack }
\newcommand{\mycomment }[1]{}

\title{On FRT--\,Clifford Algebras}
\author{I.~Heckenberger
\and A.~Sch\"uler\thanks{both authors
are supported by the Deutsche Forschungsgemeinschaft}}
\date{\small Universit\"at Leipzig, Mathematisches Institut,\\
Augustusplatz 10--11, 04109 Leipzig, Germany\\
{\footnotesize e-mail:heckenbe,schueler@mathematik.uni-leipzig.de}}

\begin{document}

\maketitle

\begin{abstract}
We study the $q$-Clifford algebras $\Cl $, called FRT--\,Clifford algebras,
introduced by Faddeev, Reshetikhin and Takhtajan.
It is shown that $\Cl $ acts on the $q$-exterior algebra $\LOq N$.
Moreover, explicit formulas for the embedding of $\Uqso N$ into
$\Cl $ and its relation to the vector and spin representations of
$\Uqso N$ are given and proved.


\end{abstract}

\section{Introduction}

The theory of spinors gives one of the most powerful structures
in differential geometry and theoretical physics. Its origin can be traced
back to works of Cartan, Dirac, Brauer and Weyl.
On the other hand, quantum groups and quantum spaces appeared
in the eighties as quantized algebras of functions on complex manifolds.
Because of their close relation to Lie groups and homogeneous spaces
one can ask whether one can develop a spin geometry of quantum spaces.
In the last years there were made several attempts in this direction.
Quantum Clifford algebras are introduced and studied for instance
in \cite{a-Hayashi90}, \cite{inp-BPR93}, \cite{a-Hannabuss00},
\cite{a-AblFau99}, \cite{a-Fiore00}.

In \cite{a-Woro2} Woronowicz presented a theory of covariant differential
calculus over Hopf algebras.
Using this or related theories quantum Clifford algebras and
spinors are defined and investigated by many authors, see
\cite{a-FadResTak1}, \cite{inp-DingFrenkel94}, \cite{a-DurdOzi94},
\cite{a-BCDRV96}, \cite{a-Heck00b}.
The approach of Woronowicz causes two important facts which are absent
in the classical situation. First, in general there exists no canonical
differential calculus on a given quantum space. Given a differential calculus,
the corresponding external algebra of differential forms depends
extremely on the algebra of functions and on the first order calculus.
Secondly, because of the noncommutative structure defining properties of
a metric tensor are much more restrictive than in the commutative situation.
Hence the varieties of Clifford algebras and spinors in the classical and
quantum setting are very different.

Exterior algebras $\LOq N$ of the quantum Euclidean space $\Oq N$
are introduced in \cite{a-FadResTak1}, \cite{a-Takeuchi90} and
\cite{a-Hayashi92}. They naturally appeared in
\cite{a-CSchW91} and \cite{a-Welk00} where the higher order calculi of the
coordinate algebra
of the quantum Euclidean space
and the quantum Euclidean sphere, respectively, were examined. The
corresponding quantum Clifford algebras, called FRT--\,Clifford algebras,
are defined in \cite{a-FadResTak1} and analyzed in \cite{inp-DingFrenkel94}.
The aim of this paper is to prove further
results on FRT--\,Clifford algebras and to present elementary methods
for the dealing with these objects.
We construct spin representations of $\Uqso N$ using FRT--\,Clifford
algebras. Related quantum spin bundles on quantum spaces
will be examined in a forthcoming paper.

In Section \ref{sec-FRTCl} the definition of the FRT--\,Clifford algebra
$\Cl $ is recalled. In Theorem \ref{t-clrep} we prove that
the algebra $\Cl $ has a representation over the
quantum exterior algebra $\LOq N$. Elementary proofs of the
semisimplicity of the algebras $\Cl $ are given.
Canonical minimal left- and right-ideals
of $\Cl $ are introduced. In Section \ref{sec-spinrep} we formulate and
give the proof of two theorems. In Theorem \ref{t-homuqgcl}
an algebra map $\pi $ from $\Uqso N$ to $\Cl $ is given.
In Theorem \ref{t-uqgreps} invariance of the
vector space of generators of $\Cl $ under the adjoint action of
$\pi (\Uqso N)$ is proved. Further, the spin representations of $\Uqso N$
are realized on the canonical left ideal of $\Cl $.

For the set of nonzero real and complex numbers we take the symbols
$\realx $ and $\compx $, respectively. Throughout we use
Einsteins convention to sum over repeated indices.
We write $\delta ^i_j=\delta _{ij}$ for Kronecker's symbol.
For integer numbers $j$, $1\leq j\leq N$, the symbol $j'$ means the number
$N+1-j$. For $m,p\in \mathbb{N}_0$ we set
\begin{align}
\qnum{m}&=\frac{1-q^{-4m}}{1-q^{-4}},&
\qnum{m}!&=\prod _{k=1}^m\qnum{k},& \qnum{0}!&=1,&
\bigqbinom{m}{p}&=\frac{\qnum{m}!}{\qnum{p}!\qnum{m-p}!}.
\end{align}
We use the definition of
the Birman--Wenzl--Murakami algebra $\mathrm{BWM}(q^{N-1},q)_k$
with $2k-2$ generators $g_i$, $e_i$, $i=1,2,\ldots ,k-1$,
given in \cite{a-BirWen89}, in the notation of \cite{a-HeckSchu99}.
The vector space of intertwiners of corepresentations $v,v'$ of the
coordinate Hopf algebra $\OOq N$ of the quantum group $\Oq N$
is denoted by $\Mor (v,v')$.
Let $u$ denote the fundamental (vector) corepresentation of $\OOq N$.
For the generators of $\Mor (u\ot u,u\ot u)$, $\Mor (1,u\ot u)$ and
$\Mor (u\ot u,1)$ we use the following symbols in the graphical calculus
(see \cite{a-CSchW91}):
\begin{align*}
\begin{xy}
0;/r1.2pc/: *{\hat{R}=}!R="mitte" +(.5,.5) \vcross ,
"mitte"+(3,0) *!L{\hat{R}^{-1}=}!R="mitte" +(.5,.5) \vcrossneg ,
"mitte"+(3,0) *!L{K=}!R="mitte" +(.5,.5) \vuncross ,
"mitte"+(3,0) *!L{C^{ij}=}!R="mitte" +(.5,.5) \vcap- ,
"mitte"+(3,0) *!L{C_{ij}=}!R="mitte" +(.5,-.5) \vcap
\end{xy}
\end{align*}

\section{The FRT--\,Clifford algebras}
\label{sec-FRTCl}

Let $q$ be a nonzero complex number, $q^m\not=1 $ for all $m\in \mathbb{N}$.
We set $\qp =q+q^{-1}$ and
$\qm =q-q^{-1}$. Let $N\in \mathbb{N}$, $N\geq 3$, and $n\in \mathbb{N}$ such
that $N=2n+\epsilon $, $\epsilon \in \{0,1\}$.
In the sequel we frequently use the following standard notation
for $\Rda $-matrices corresponding to orthogonal quantum groups
$\Oq N$, cf.~\cite[Subsect.~1.4]{a-FadResTak1}.
\begin{equation}
\begin{gathered}
\Rda {}^{ij}_{kl}=q^{2\delta ^i_j-2\delta ^i_{j'}}\delta ^i_l\delta ^j_k
+(q^2-q^{-2})\delta _{i<l}(\delta ^i_k\delta ^j_l-K^{ij}_{kl}),\\
C^{ij}=C_{ij}=q^{-2\rho _i}\delta ^i_{j'},\quad K^{ij}_{kl}=C^{ij}C_{kl},
\qquad \rho _i=\delta _{i<i'}(N/2-i)-\delta _{i'<i}(N/2-i'),\\
P_+=\frac{1}{q^2+q^{-2}}\left(q^{-2}\,\id +\Rda -\frac{\qp \qm }{q^{2N}-1}
K\right),\\
P_-=\frac{1}{q^2+q^{-2}}\left(q^2\,\id -\Rda -\frac{\qp \qm }{q^{2N-4}+1}
K\right),
\end{gathered}
\end{equation}
where $\delta _{i<l}=1$ if $i<l$ and $\delta _{i<l}=0$ otherwise.

Let $c^2\in \realx $.
By Definition 13 in \cite{a-FadResTak1} the quantum Clifford algebra
$\Cl $ is the complex unital algebra generated by the elements
$\gcl _i$, $i=1,\ldots ,N$, and relations
\begin{align}\label{eq-defclrel}
P_+{}^{kl}_{ij}\gcl _k\gcl _l&=0,\quad
i,j=1,2,\ldots ,N,&  C^{ij}\gcl _i\gcl _j&=
c^2\frac{q^{2N}-1}{q^2-1}.
\end{align}
The factor in the last relation is choosen in such a way that
formulas in the sequel become more simple.
The explicit form of the relations of $\Cl $ can be computed
similarly to the proof of Proposition 9.14 in \cite{b-KS}. We obtain
\begin{equation}\label{eq-cliff}
\begin{gathered}
\gcl _i\gcl _i=0,\quad i\not= i',\qquad
\gcl _j\gcl _i=-q^2\gcl _i\gcl _j,\quad j>i, j\not= i',\\
\gcl _{i'}\gcl _i=-\gcl _i\gcl _{i'}+(q^2-q^{-2})\sum _{j=1}^{i-1}q^{2j-2i+2}
\gcl _j\gcl _{j'}+c^2q^{N-2i+1}\qp ,\quad i<i',\\
\epsilon \gcl _{n+1}\gcl _{n+1}=\epsilon \left(\qm \sum _{j=1}^n
q^{2j-2n}\gcl _j\gcl _{j'}+c^2\right).
\end{gathered}
\end{equation}
It is easy to see that the mapping $\gcl _i\mapsto c\gcl _i$ defines an
algebra isomorphism $\Cl \to \Clg{1}$.
In contrast, the choice $c=0$ leads to the $q$-exterior algebra $\LOq N$
which has been studied for instance in \cite{a-FadResTak1},
\cite{a-Takeuchi90} or \cite[Sect. 5]{a-Hayashi92}.

The algebra $\Cl $ obeys a $\mathbb{Z}_2$-grading $\partial _0$
given by $\partial _0(\gcl _i)=1$ for all $i=1,\ldots ,N$. The components
of even and odd degree are denoted by $\Cle $ and $\Clo $, respectively.
There exists also a $\mathbb{Z}^n$-grading $\partial =(\partial _1,
\partial _2,\ldots ,\partial _n)$ of $\Cl $ which is determined
by $\partial _i(\gcl _j)=\delta ^i_j-\delta ^i_{j'}$ for $i=1,\ldots ,n$
and $j=1,\ldots ,N$. For $N=2n+1$ the mapping $\partial _{n+1}:
\gcl _j\mapsto \delta ^j_{n+1}$, $j=1,2,\ldots ,N$, defines another
$\mathbb{Z}_2$-grading of $\Cl[2n+1]$.

Let $\alpha _i$, $i=1,\ldots ,N$, be nonzero complex numbers such that
$\alpha _i\alpha _{i'}=1$. From (\ref{eq-cliff}) it
follows that the mapping $\gcl _i\mapsto \alpha _i\gcl _i$
defines an algebra isomorphism of $\Cl $.
Moreover, there exists an algebra antiautomorphism $\tau $ of
$\Cl $ given by $\tau (\gcl _i)=\gcl _{i'}$. This
antiautomorphism as well as its compose with any of the above automorphisms
has order two.

Let $V$ denote the complex vector space with basis
$\{\gcl _i\,|\,i=1,\ldots ,N\}$ and let $V^{\ot k}$, $k\geq 1$, be the
$k$-fold tensor product $V\ot V\ot \cdots \ot V$.
We write $V^\ot $ for $\bigoplus _{k=0}^\infty V^{\ot k}$,
where $V^{\ot 0}=\comp 1$.
The mapping $\varphi :\mathrm{BWM}(q^{2N-2},q^2)_{k+1}\to \End{V^{\ot k+1}}$,
given by
$\varphi (g_i)=\Rda _{i,i+1}$, $\varphi (e_i)=K_{i,i+1}$, $i=1,\ldots ,k$,
is a representation of $\mathrm{BWM}(q^{2N-2},q^2)_{k+1}$. If no confusion
can arise we will identify elements of the algebra
$\mathrm{BWM}(q^{2N-2},q^2)_{k+1}$
with their images under the representation $\varphi $.

In \cite{a-HeckSchu99}, formulas (19) and (20), the elements $d'{}^-_{k+1,i}$
and $b^-_{1,k}$, $k=0,1,\ldots ,N-1$, $i=1,\ldots ,k$, of the
Birman--Wenzl--Murakami algebra $\mathrm{BWM}(q^{2N-2},q^2)_{k+1}$ were
introduced.
\begin{equation}\label{eq-bm}
\begin{aligned}
d'{}^-_{k+1,i}&=K_{12}K_{23}\cdots K_{k-i+1,k-i+2}
\sum _{j=0}^{i-1}(-q^{-2})^j\prod _{l=1}^j\Rda _{k-i+1+l,k-i+2+l},\\
b^-_{1,k}&=\sum _{i=0}^k(-q^{-2})^i\Rda _{12}\Rda _{23}\cdots \Rda _{i,i+1}
-\frac{\qp \qm }{1+q^{2N-4k}}\sum _{i=1}^k q^{4i-4k-2}d'{}^-_{k+1,i}.
\end{aligned}
\end{equation}
Using (\ref{eq-bm}) the unique antisymmetrizer $A_{k+1}$ of
$\mathrm{BWM}(q^{2N-2},q^2)_{k+1}$ is constructed.
Recall that $A_1=\id $, $A_2=P_-$ and
$A_{i+1}=(\id \ot A_i)b^-_{1,i}/\qnum{i+1}$ for $i>2$.
Further, from $A_{m+1}=A_{m+1}(\id \ot A_m)$,
$\Rda _{i,i+1}A_{m+1}=-q^{-2}A_{m+1}$,
$K_{i,i+1}A_{m+1}=0$ ($i=1,2,\ldots ,m$) and (\ref{eq-bm}) we conclude that
\begin{equation}\label{eq-antsymrec}
\begin{aligned}
A_{m+1}&=\frac{1}{\qnum{m+1}}(\id \ot A_m)b^-_{1,m}=\frac{1}{\qnum{m+1}}
(\id \ot A_m)b^-_{1,m}(\id \ot A_m)\\
&=\frac{1}{\qnum{m+1}}(\id \ot A_m)\left(\id -q^{-2}\qnum{m}\Rda _{12}
-\frac{1-q^{-4m}}{1+q^{2N-4m}}K_{12}\right) (\id \ot A_m).
\end{aligned}
\end{equation}
We use the symbol $V^{\wedge k}$ for the vector space $V^{\ot k}/\ker A_k$.
Then $\LOq N=\bigoplus _{k=0}^NV^{\wedge k}$ as vector spaces and the
multiplication $\mwedge $ corresponds to the tensor product in $V^\ot $.

Let $g:V\ot V\to \comp $ denote the linear mapping given by the formula
$g(\gcl _i\ot \gcl _j)=c^2q^{2N-3}\qp C_{ij}/(q^{2N-4}+1)$.
{}From $C^{ij}C_{ij}=q^{-2N+2}(q^{2N}{-}1)(q^{2N-4}{+}1)/(q^2{-}q^{-2})$ we
then get
$a_{ij}\gcl _i\gcl _j=g(a_{ij}\gcl _i\ot \gcl _j)$ in $\Cl $
whenever $a_{ij}\gcl _i\mwedge \gcl _j=0$ in $\LOq N$.
Inspired by the methods presented in \cite{a-Heck99} and \cite{a-BCDRV96},
let us define contraction mappings $\ctr{\cdot }{\cdot }:V^{\ot k}\ot
V^{\ot l}\to V^{\ot \mathrm{max}\{0,l-k\}}$ by $\ctr{\rho _k}{\rho '_l}=0$ for
$k>l$ and
\begin{align}
\ctr{\rho }{\rho '_l}&=g_{12}(\rho \ot b^-_{1,l-1}(\rho '_l)),&
\ctr{\rho \ot \rho _k}{\rho '_l}&=\ctr{\rho }{\ctr{\rho _k}{\rho '_l}}
\end{align}
for $\rho \in V$, $\rho _i,\rho '_i\in V^{\ot i}$.
Now let $l\geq k$, $\rho '_l\in V^{\ot l}$,
$\rho _k:=a_J\gcl _{j_1}\gcl _{j_2}\cdots \gcl _{j_k}\in V^{\ot k}$,
$a_J\in \comp $ for $J=(j_1,\ldots ,j_k)$, $1\leq j_m\leq N$ for all $m$.
We show that $A_{l-k}\ctr{\rho _k}{\rho '_l}=0$ provided $A_l(\rho '_l)=0$ or
$A_k(\rho _k)=0$.
Indeed, by the definition of the contraction mapping we get
\begin{align*}
A_{l-k}\ctr{\rho _k}{\rho '_l}&=a_JA_{l-k}\ctr{\gcl _{j_1}}{\ctr{%
\gcl _{j_2}}{\cdots \ctr{\gcl _{j_k}}{\rho '_l}}}\\
&=a_Jg_{12}(\gcl _{j_1}\ot \qnum{l-k+1}A_{l-k+1}(\ctr{\gcl _{j_2}}{\cdots
\ctr{\gcl _{j_k}}{\rho '_l}}))=\cdots \\
&=\frac{\qnum{l}!}{\qnum{l-k}!}a_Jg_{12}(\gcl _{j_1}\ot g_{12}(\gcl _{j_2}\ot 
\cdots g_{12}(\gcl _{j_k}\ot A_l(\rho '_l)))).
\end{align*}
Hence $\ctr{\cdot }{\cdot }$ maps $V^{\ot k}\ot V^{\wedge l}$, $l\geq k$, to
$V^{\wedge l-k}$.
On the other hand, the last line of the above equation can also be written in
the form
\begin{align*}
A_{l-k}\ctr{\rho _k}{\rho '_l}&=\frac{\qnum{l}!}{\qnum{l-k}!}
g_{12}g_{23}\cdots g_{k,k+1}(\rho _k\ot A_l(\rho '_l)).
\end{align*}
Suppose now that $A_k \rho _k=0$. Then $\rho _k$ is a linear combination of
elements $\xi _0:=\xi '_i\ot \xi \ot \xi ''_{k-i-2}$, $\xi '_i\in V^{\ot i}$,
$\xi ''_{k-i-2}\in V^{\ot k-i-2}$,
$\xi \in V\ot V$, $A_2(\xi )=0$. Without loss of generality
we can assume that $\Rda \xi =r\xi $ with $r=q^2$ or $r=q^{-2N+2}$. This means
that $\xi $ is an element of the image of the projector $P_+$ or $P_0$,
respectively.  For such elements we easily obtain
\begin{align*}
A_{l-k}\ctr{\xi _0}{\rho '_l}&=A_{l-k}\ctr{r^{-1}\Rda _{i+1,i+2}\xi _0}{%
\rho '_l}\\
&=\frac{r^{-1}\qnum{l}!}{\qnum{l-k}!}
g_{12}g_{23}\cdots g_{k,k+1}(\Rda _{i+1,i+2}\xi _0\ot A_l(\rho '_l))\\
&=\frac{r^{-1}\qnum{l}!}{\qnum{l-k}!}
g_{12}g_{23}\cdots g_{k,k+1}(\xi _0\ot \Rda _{k-i-1,k-i}A_l(\rho '_l))\\
&=\frac{-q^{-2}r^{-1}\qnum{l}!}{\qnum{l-k}!}
g_{12}g_{23}\cdots g_{k,k+1}(\xi _0\ot A_l(\rho '_l))
=-q^{-2}r^{-1}A_{l-k}\ctr{\xi _0}{\rho '_l}.
\end{align*}
Since $q$ is not a root of unity, we obtain
$A_{l-k}\ctr{\xi _0}{\rho '_l}=0$
and hence $A_{l-k}\ctr{\rho _k}{\rho '_l}=0$.
This means that $\ctr{\cdot }{\cdot }$ is a well-defined mapping from
$V^{\wedge k}\ot V^{\wedge l}$, $k\leq l$, to $V^{\wedge l-k}$.

\begin{thm}\label{t-clrep}
There exists a representation $\lact $ of $\Cl $ on $\LOq N$
such that
\begin{align}\label{eq-clrep}
\gcl _i\lact \rho _m &=\gcl _i\mwedge \rho _m+\frac{1+q^{2N-4}}{1+q^{2N-4m}}
\ctr{\gcl _i}{\rho _m}
\end{align}
for $\rho _m\in V^{\wedge m}$, $m\in \{0,1,\ldots ,N\}$, $i=1,2,\ldots ,N$.
\end{thm}

\begin{bew}
Since $\ctr{\cdot }{\cdot }$ is a well-defined mapping from
$V\ot V^{\wedge m}$ to
$V^{\wedge m-1}$, we only have to show that the definition (\ref{eq-clrep})
is compatible with the relations (\ref{eq-defclrel}) of $\Cl $.
Suppose that $m=0$. Then
\begin{align*}
\gcl _i\gcl _j\lact 1&=\gcl _i\lact \gcl _j=\gcl _i\mwedge \gcl _j
+\ctr{\gcl _i}{\gcl _j}
\end{align*}
for $i,j=1,\ldots ,N$. Since $g(P_+{}_{ij}^{kl}\gcl _k\ot \gcl _l)=0$
we obtain $P_+{}_{ij}^{kl}\gcl _k\gcl _l\lact 1=0$.
{}From $C^{ij}C_{ij}=q^{-2N+2}(q^{2N}-1)(q^{2N-4}+1)/(q^2-q^{-2})$ we
conclude that
\begin{align*}
C^{ij}\gcl _i\gcl _j\lact 1&=C^{ij}g(\gcl _i\ot \gcl _j)=
C^{ij}\frac{c^2q^{2N-3}\qp }{q^{2N-4}+1}C_{ij}=\frac{c^2q^{-1}(q^{2N}-1)}{\qm }
\end{align*}
which proves the compatibility with the second equation of (\ref{eq-defclrel})
for $m=0$.

Now consider the case $m\geq 1$. Then we get
\begin{align*}
a_{ij}\gcl _i\gcl _j\lact \rho _m&=a_{ij}\gcl _i\lact \left(\gcl _j\mwedge
\rho _m
+\frac{1+q^{2N-4}}{1+q^{2N-4m}}\ctr{\gcl _j}{\rho _m}\right)\\
&=a_{ij}\gcl _i\mwedge \gcl _j\mwedge \rho _m
+\frac{(1+q^{2N-4})^2}{(1+q^{2N-4m})(1+q^{2N-4m+4})}\ctr{a_{ij}\gcl _i\ot
\gcl _j}{\rho _m}\\
&\phantom{=}
+\frac{1+q^{2N-4}}{1+q^{2N-4m}}a_{ij}\gcl _i\mwedge \ctr{\gcl _j}{\rho _m}
+\frac{1+q^{2N-4}}{1+q^{2N-4m-4}}a_{ij}\ctr{\gcl _i}{\gcl _j\ot \rho _m}
\end{align*}
for $\rho _m\in V^{\wedge m}$, $i,j=1,\ldots ,N$, $a_{ij}\in \comp $.
If $a_{ij}\gcl _i\ot \gcl _j\in \ker A_2$ then the first two summands of the
last equation vanish and we get an element in $V^{\wedge m}$. More precisely,
for $\rho _m=\gcl _{k_1}\mwedge \gcl _{k_2}\mwedge \cdots \gcl _{k_m}$
the expression $A_m(a_{ij}\gcl _i\gcl _j\lact \rho _m)$ takes the form
\begin{gather*}
\epsfig{file=frtclif.bilder.001,width=\linewidth }
\end{gather*}
Using $A_m=A_m(\id \ot A_{m-1})$ and $(\id \ot A_{m-1})b^-_{1,m-1}=
\qnum{m}A_m$ in the first summand and $(\id \ot A_m)b^-_{1,m}=
\qnum{m+1}A_{m+1}$ in the second one the coefficient in the above
expression becomes
\begin{gather*}\tag{$*$}
\epsfig{file=frtclif.bilder.002}
\end{gather*}
Now in the second summand we replace $A_{m+1}$ by the second line of
(\ref{eq-antsymrec}). Then for $a_{ij}=C^{ij}$ the expression ($*$)
becomes
\begin{align*}
\biggl(\frac{c^2q^{2N-3}\qp \qnum{m}}{q^{2N-4m}+1}&+
\frac{c^2q^{2N-3}\qp }{q^{2N-4m-4}+1}\biggl(
\frac{q^{-2N+2}(q^{2N}-1)(q^{2N-4}+1)}{\qp \qm }-\\
&-q^{-2}\qnum{m}q^{2N-2}-\frac{1-q^{-4m}}{1+q^{2N-4m}}\biggr)
\biggr)A_m=\frac{c^2(q^{2N}-1)}{q^2-1}A_m
\end{align*}
which proves the compatibility of $\lact $ with the last relation in
(\ref{eq-defclrel}).
Setting $a_{ij}=P_+{}_{rs}^{ij}$, $r,s=1,\ldots ,N$, ($*$) results in
\begin{align*}
\biggl(\frac{c^2q^{2N-3}\qp \qnum{m}}{q^{2N-4m}+1}&+
\frac{c^2q^{2N-3}\qp }{q^{2N-4m-4}+1}\biggl(
0-q^{-2}\qnum{m}q^{-2}-\frac{1-q^{-4m}}{1+q^{2N-4m}}\biggr)
\biggr)M_{+,m}=0,
\end{align*}
where
\begin{gather*}
\epsfig{file=frtclif.bilder.003}
\end{gather*}
\end{bew}

Let $I_N=(i_1,i_2,\ldots ,i_N)$ and $J_n=(j_1,j_2,\ldots ,j_n)$,
$i_k,j_k\in \{0,1\}$ for all $k$.
We use the abbreviation
$\gcl ^{I_N}=\gcl _1^{i_1}\gcl _2^{i_2}\cdots \gcl _N^{i_N}$.
Set $\ell (I_N)=\sum _{k=1}^Ni_k$,
$\ell (J_n)=\sum _{k=1}^nj_k$, and $L(J_n)=\sum _{k=1}^nkj_k$.

\begin{satz}\label{s-clbasis}
The set $\{\gcl _1^{i_1}\gcl _2^{i_2}\cdots \gcl _N^{i_N}\,|\,
i_k\in \{0,1\}, k=1,\ldots ,N\}$ forms a vector space basis of
$\Cl $.
\end{satz}

\begin{bew}
{}From the explicit form (\ref{eq-cliff}) of the generating relations of
the algebra $\Cl $ we conclude that the elements $\gcl ^{I_N}$
span the vector space $\Cl $. To prove that they form a basis
we use Theorem \ref{t-clrep}. Suppose that there is an element
$\rho =\sum _{I_N}\lambda _{I_N}\gcl ^{I_N}=0$ in $\Cl $.
Then $\rho \lact 1=0$. Let
$j_{\mathrm{max}}=\mathrm{max}\{\ell (I_N)\,|\,\lambda _{I_N}\not=0\}$.
Then we have
\begin{align*}
\rho \lact 1=\sum _{J_N,\ell (J_N)=j_{\mathrm{max}}}\lambda _{J_N}
\gcl _1^{j_1}\mwedge \gcl _2^{j_2}\mwedge \cdots \mwedge \gcl _N^{j_N}+\rho ',
\end{align*}
where $\rho '\in \bigoplus _{k=0}^{j_{\mathrm{max}}-1}V^{\wedge k}$.
Since $\rho \lact 1=0$, we conclude that $\lambda _{J_N}=0$ for all $J_N$ with
$\ell (J_N)=j_{\mathrm{max}}$. This is a contradiction to the definition of
$j_{\mathrm{max}}$.
\end{bew}

\begin{folg}
The representation $\lact $ of $\Cl $ on $\LOq N$ is faithful.
\end{folg}

Let $\Clz $ denote the subalgebra of $\Cl $ defined by
\begin{align}
\Clz =\{z\in \Cl \,|\,\partial _k(z)=0 \quad
\text{for all $k=1,2,\ldots ,n$}\}.
\end{align}

\begin{satz}\label{s-clzcomm}
The algebra $\Clz $ is commutative.
\end{satz}

\begin{bew}
By Proposition \ref{s-clbasis} the set
\begin{align}
\{\gcl _1^{i_1}\gcl _2^{i_2}\cdots \gcl _n^{i_n}\gcl _{n+1}^{i_{n+1}}
\gcl _{n'}^{i_n}\cdots \gcl _{N-1}^{i_2}\gcl _N^{i_1}\,|\,
i_l\in \{0,1\},l=1,2,\ldots ,n+1,\,i_{n+1}\leq \epsilon \}
\end{align}
forms a vector space basis of the algebra $\Clz[2n+\epsilon ]$.
{}From this fact and the first line in (\ref{eq-cliff}) we conclude that
$\Clz $ is generated as an algebra by the set of its elements
\begin{align*}
&\{\gcl _i\gcl _{i'}\,|\,i=1,2,\ldots ,n\}& &\text{for $N=2n$,}\\
&\{\gcl _{n+1},\gcl _i\gcl _{i'}\,|\,i=1,2,\ldots ,n\}\quad &
&\text{for $N=2n+1$.}
\end{align*}
Now the commutativity of these generators follows immediately from
the equations in the first line of (\ref{eq-cliff}).
\end{bew}

The following proposition, as well as Proposition \ref{s-oddsemis}
was proved in \cite{inp-DingFrenkel94} using representation theory of
$\Uqso N$. Here we give an elementary proof.

\begin{satz}\label{s-evensimple}
The algebra $\Cl[2n]$ is simple.
\end{satz}

\begin{bew}
Let $\mc{I}\not=\{0\}$ be a two-sided ideal of $\Cl[2n]$. We prove that
$\mc{I}=\Cl[2n]$. First it will be shown by induction that
\begin{itemize}
\item[($*$)]
for all $j\in \{1,2,\ldots ,n\}$ there exists
$\rho _j=\sum _{I_N}\lambda _{I_N}\gcl ^{I_N}\in \mc{I}$ such that
$\lambda _{I_N}\not=0$ implies $i_k=1$ for all $k<j$.
\end{itemize}
For $j=1$ this only means
that there exists a nonzero $\rho _1\in \mc{I}$. Further, from
(\ref{eq-cliff}) we obtain
\begin{align*}
\gcl _j\rho _j&=\sum _{I_N,i_j=0}(-q^2)^{j-1}\lambda _{I_N}\gcl _1\gcl _2
\cdots \gcl _j \gcl _{j+1}^{i_{j+1}}\cdots \gcl _N^{i_N}.
\end{align*}
If $\gcl _j\rho _j=0$ then Proposition \ref{s-clbasis} gives
$\lambda _{I_N}=0$ for all $I_N$ with $i_j=0$.
Hence $\rho _{j+1}:=\rho _j$ fulfills our hypothesis
($*$) for $j+1$.
Otherwise $\gcl _j\rho _j\not=0$ and the element $\rho _{j+1}:=\gcl _j\rho _j$
satisfies ($*$) for $j+1$. By induction it follows that ($*$) is valid for
$j=n+1$.

Now we prove by induction that
\begin{itemize}
\item[($**$)]
for all $j\in \{1,2,\ldots ,n\}$ there exists
$\rho _j=\sum _{I_N}\lambda _{I_N}\gcl ^{I_N}\in \mc{I}$ such that
$\lambda _{I_N}\not=0$ implies $i_k=1$ for $k\leq n$ and for $k>j'$.
\end{itemize}
By ($*$) this is true for $j=1$. Further, we obtain
\begin{align*}
\rho _j\gcl _{j'}&=\sum _{I_N,i_{j'}=0}(-q^2)^{j-1}\lambda _{I_N}\gcl _1^{i_1}
\cdots \gcl _{j'-1}^{i_{j'-1}}\gcl _{j'}\gcl _{j'+1}\cdots \gcl _N.
\end{align*}
If $\rho _j\gcl _{j'}=0$ then Proposition \ref{s-clbasis} implies
$\lambda _{I_N}=0$ whenever
$i_{j'}=0$. Hence $\rho _{j+1}:=\rho _j$ satisfies ($**$). Otherwise
$\rho _{j+1}:=\rho _j\gcl _{j'}\not=0$
and it fulfills ($**$). By induction we get
($**$) for $j=n+1$. This means that $\gcl _1\gcl _2\cdots \gcl _N\in \mc{I}$.

Now we prove by induction that $\gcl _1\gcl _2\cdots \gcl _{j'}\in \mc{I}$
for $j\leq n+1$. Indeed, this is true for $j=1$. Further, from
(\ref{eq-cliff}) and $j\leq n$ we get
\begin{align*}
\gcl _1\gcl _2\cdots \gcl _{j'}\gcl _j&=
\gcl _1\gcl _2\cdots \gcl _{j'-1}\left(-\gcl _j\gcl _{j'}+\qp \qm
\sum _{k=1}^{j-1}q^{2k-2j+2}\gcl _k\gcl _{k'}+c^2q^{N-2j+1}\qp \right)\\
&=c^2q^{N-2j+1}\qp \gcl _1\gcl _2\cdots \gcl _{j'-1}.
\end{align*}
This means that $\gcl _1\gcl _2\cdots \gcl _n\in \mc{I}$. Finally, we prove
by inverse induction on $j$
that $\rho _j:=\gcl _1\gcl _2\cdots \gcl _j\in \mc{I}$ for $j=n,n-1,\ldots ,0$.
This is true for $j=n$.
Of course, then $\rho _j\gcl _{j'}\in \mc{I}$. On the other hand,
\begin{align*}
\gcl _{j'}\rho _j&=
(-q^2)^{j-1}\gcl _1\gcl _2\cdots \gcl _{j-1}\gcl _{j'}\gcl _j\\
&=(-q^2)^{j-1}\gcl _1\cdots \gcl _{j-1}\left(-\gcl _j\gcl _{j'}+
\qp \qm \sum _{k=1}^{j-1}q^{2k-2j+2}\gcl _k\gcl _{k'}+c^2q^{N-2j+1}\qp
\right)\\
&=-(-q^2)^{j-1}\rho _j\gcl _{j'}+(-1)^{j-1}c^2q^{N-1}\qp
\gcl _1\gcl _2\cdots \gcl _{j-1}.
\end{align*}
Hence $\rho _{j-1}\in \mc{I}$. For $j=0$ we get $\rho _0=1\in \mc{I}$
and hence $\mc{I}=\Cl $.
\end{bew}

For $N=2n+1$ let $\mc{I}_+$ and $\mc{I}_-$ denote the two-sided ideal of
$\Cl $ generated by the element $\rho _+$ and $\rho _-$,
respectively, where
\begin{align}
\rho _\pm =\gcl _1\cdots \gcl _n(\pm c+\gcl _{n+1})\gcl _{n+2}\cdots \gcl _N.
\end{align}

\begin{satz}\label{s-oddsemis}
The algebra $\Cl[2n+1]$ is semisimple and
splits into the direct sum of the two simple ideals
$\mc{I}_+$ and $\mc{I}_-$. Further, $\dim \mc{I}_+=\dim \mc{I}_-$.
\end{satz}

\begin{bew}
Let $\mc{I}\not=\{0\}$ be a two-sided ideal of $\Cl[2n+1]$.
Similarly to the first part of the proof of Proposition \ref{s-evensimple}
one can show that there exist complex numbers $\alpha _1,\alpha _2$,
$|\alpha _1|+|\alpha _2|>0$, such that
$\rho :=\gcl _1\gcl _2\cdots \gcl _n(\alpha _1+\alpha _2\gcl _{n+1})
\gcl _{n+2}\cdots \gcl _N \in \mc{I}$. Multiplying $\rho $ by
$(-q^{-2})^n\gcl _{n+1}$ from the left we obtain that
\begin{align*}
\rho ':=\gcl _1\cdots \gcl _n\gcl _{n+1}(\alpha _1+\alpha _2\gcl _{n+1})
\gcl _{n+2}\cdots \gcl _N=
\gcl _1\cdots \gcl _n(c^2\alpha _2+\alpha _1\gcl _{n+1})
\gcl _{n+2}\cdots \gcl _N
\end{align*}
belongs to $\mc{I}$. If $\rho $ and $\rho '$ are linearly independent then
$\mc{I}$ contains both $\mc{I}_+$ and $\mc{I}_-$. Otherwise
$\alpha _1^2=c^2\alpha _2^2$ and $\mc{I}$ contains at least one of the
ideals $\mc{I}_+$ and $\mc{I}_-$. Assume $\mc{I}_\eta \subset \mc{I}$
with $\eta =+$ or $\eta =-$.  Since $\Cl[2n+1]$
is a finite dimensional algebra over $\comp $, its radical is nilpotent.
Suppose that $\mc{I}_\eta $ is contained in the radical of $\Cl[2n+1]$.
Then $\rho _\eta $ has to be nilpotent. Let us determine $\rho _\eta ^2$.
First, for $j=1,\ldots ,n$ we compute
\begin{align*}
&\gcl _1\cdots \gcl _n(\eta c+\gcl _{n+1})\gcl _{n+2}\cdots \gcl _{j'}\cdot
\gcl _j\gcl _{j+1}\cdots \gcl _n(\eta c+\gcl _{n+1})=\\
&\qquad
=\gcl _1\cdots \gcl _n(\eta c+\gcl _{n+1})\gcl _{n+2}\cdots \gcl _{j'-1}
\biggl(-\gcl _j\gcl _j'+\qp \qm \sum _{k=1}^{j-1}q^{2k-2j+2}\gcl _k\gcl _{k'}\\
&\qquad \phantom{=}+c^2q^{N-2j+1}\qp \biggr)\gcl _{j+1}\cdots \gcl _n
(\eta c+\gcl _{n+1})\\
&\qquad =c^2q^{N-2j+1}\qp \gcl _1\cdots \gcl _n(\eta c+\gcl _{n+1})\gcl _{n+2}
\cdots \gcl _{j'-1}\cdot \gcl _{j+1}\cdots \gcl _n(\eta c+\gcl _{n+1}).
\end{align*}
Further,
\begin{align*}
\gcl _1\cdots \gcl _n(\eta c+\gcl _{n+1})^2&=
\gcl _1\cdots \gcl _n\biggl(c^2+2\eta c\gcl _{n+1}+\qm \sum _{k=1}^nq^{2k-2n}
\gcl _k\gcl _{k'}+c^2\biggr)\\
&=
2\eta c\gcl _1\cdots \gcl _n(\eta c+\gcl _{n+1}).
\end{align*}
This yields
\begin{align}\label{eq-rhogamma}
\rho _\eta \gcl _1\gcl _2\cdots \gcl _n(\eta c+\gcl _{n+1})&=
2\eta c^{2n+1}q^{n(n+1)}\qp ^n\gcl _1\gcl _2\cdots \gcl _n(\eta c+\gcl _{n+1})
\end{align}
and hence $\rho _\eta^2=2\eta c^{2n+1}q^{n(n+1)}\qp ^n\rho _\eta $.
Therefore $\rho _\eta $ is not
nilpotent. The only possibility is that the radical is $\{0\}$ and hence
$\Cl $ is semisimple.

Since $\dim \Cl[2n+1]=2^{2n+1}\not= k^2$ for all $k\in \mathbb{N}$,
$\Cl[2n+1]$ is not simple and $\mc{I}_\eta \not=\Cl[2n+1]$ for $\eta =+$
or $\eta =-$. Then the complement ideal $\Cl \circleddash \mc{I}_\eta $
contains $\mc{I}_{-\eta }$ and we get
$\Cl[2n+1]\circleddash (\mc{I}_+\oplus \mc{I}_-)=\{0\}$.
Hence $\Cl[2n+1]=\mc{I}_+\oplus \mc{I}_-$.
Then the fact that each nontrivial simple ideal contains $\mc{I}_+$
or $\mc{I}_-$ in turn implies that both $\mc{I}_+$ and $\mc{I}_-$ are
simple ideals.
Moreover, since the mapping $\gcl _i\mapsto \gcl _i$, $i\not= n+1$,
$\gcl _{n+1}\mapsto -\gcl _{n+1}$ defines an algebra isomorphism
of $\Cl[2n+1]$ and maps $\mc{I}_+$ to $\mc{I}_-$, both ideals have the
same dimension.
\end{bew}

For $N=2n+1$ let $\nu \in \{1,-1\}$ and for $N=2n$ let $\nu =1$.
Let $\varphi _N^\nu \in \Cl $ denote the element
$\varphi ^\nu _N=(\nu c+\gcl _{n+1})\gcl _{n+2}\gcl _{n+3}\cdots \gcl _N$
for $N=2n+1$ and $\varphi ^1_N=\gcl _{n+1}\gcl _{n+2}\cdots \gcl _N$
for $N=2n$.

\begin{satz}\label{s-idealbasis}
The left ideals $\mc{I}^\nu _\mathrm{L}:=\Cl \varphi ^\nu _N$
of $\Cl $ are minimal. For all $\nu $
the set $\{\gcl _1^{i_1}\gcl _2^{i_2}\cdots \gcl _n^{i_n}\varphi ^\nu _N\,|\,
i_j\in \{0,1\},j=1,\ldots ,n\}$ forms a vector space basis of
$\mc{I}^\nu _\mathrm{L}$.
\end{satz}

\begin{bew}
By definition the given elements belong to the left ideal
$\mc{I}^\nu _\mathrm{L}$. On the other hand, the equations (\ref{eq-cliff})
imply that $\gcl _k\varphi ^\nu _N=0$ for $k>k'$ and
$\gcl _{n+1}\varphi ^\nu _{2n+1}=\nu c\varphi ^\nu _{2n+1}$.
Hence each expression $\gcl _1^{i_1}\gcl _2^{i_2}\cdots \gcl _N^{i_N}
\varphi ^\nu _N$, $i_l\in \{0,1\}$ for $l=1,2,\ldots ,N$, is a scalar
multiple of
$\gcl _1^{i_1}\gcl _2^{i_2}\cdots \gcl _n^{i_n}\varphi ^\nu _N$.
The linear independence of the given set follows immediately from Proposition
\ref{s-clbasis}.
Finally, the dimension of the ideal $\mc{I}_\mathrm{L}^\nu $
is $2^n$. By Propositions \ref{s-evensimple} and \ref{s-oddsemis}
and from $\dim \Cl =2^N$ it follows that
minimal left ideals have dimension $2^n$. Hence $\mc{I}^\nu _\mathrm{L}$
is minimal.
\end{bew}

Similarly there exist canonical minimal right ideals $\mc{I}^\nu _\mathrm{R}$
of $\Cl $ generated by $\psi ^1_{2n}=\gcl _1\gcl _2\cdots \gcl _n$ and
$\psi ^\nu _{2n+1}=\gcl _1\gcl _2\cdots \gcl _n(\nu c+\gcl _{n+1})$,
respectively.

\section{Spin representations of $\Uqso N$}
\label{sec-spinrep}

Let $N\in \mathbb{N}$, $N\geq 3$, $N=2n+\epsilon $, where
$n\in \mathbb{N}$ and $\epsilon \in \{0,1\}$.
If not otherwise stated we use the conventions in \cite{b-KS}.
However we set $q_i=q^{d_i}$, $i=1,2,\ldots ,n$, where
$d_i=2$ for $i<n$ and $d_n=2-\epsilon $.
In \cite{b-KS} the algebra $\Uqso N$ is defined by the generators
$E_i,F_i,K_i$ and $K_i^{-1}$, $i=1,2,\ldots ,n$, and relations
\begin{gather}
\label{eq-UK}
K_iK_j=K_jK_i,\qquad K_iK_i^{-1}=K_i^{-1}K_i=1,\\
\label{eq-UKEKF}
K_iE_jK_i^{-1}=q_i^{a_{ij}}E_j,\qquad K_iF_jK_i^{-1}=q_i^{-a_{ij}}F_j,\\
\label{eq-UEF}
E_iF_j-F_jE_i=\delta _{ij}\frac{K_i-K_i^{-1}}{q_i-q_i^{-1}},\\
\label{eq-USerreE}
\sum _{r=0}^{1-a_{ij}}(-1)^rq_i^{r(n-r)}\bigqbinom{1-a_{ij}}{r}_{q_i}
E_i^{1-a_{ij}-r}E_jE_i^r=0,\quad i\not=j,\\
\label{eq-USerreF}
\sum _{r=0}^{1-a_{ij}}(-1)^rq_i^{r(n-r)}\bigqbinom{1-a_{ij}}{r}_{q_i}
F_i^{1-a_{ij}-r}F_jF_i^r=0,\quad i\not=j,
\end{gather}
$i,j=1,2,\ldots ,n$,
where $\qbinom{m}{p}_{q^2}=\qbinom{m}{p}$, and $\qbinom{m}{p}_q$ is defined
in the same way as $\qbinom{m}{p}$ but everywhere $q^2$ has to be replaced
by $q$.

By \cite[Proposition 6.5]{b-KS} there exists a unique Hopf algebra structure
on $\Uqso N$ with coproduct $\copr $, counit $\coun $ and antipode $\antip $
such that
\begin{gather}
\copr (K_i)=K_i\ot K_i,\qquad \copr (K_i^{-1})=K_i^{-1}\ot K_i^{-1},\\
\copr (E_i)=1\ot E_i+E_i\ot K_i,\qquad
\copr (F_i)=K_i^{-1}\ot F_i+F_i\ot 1,\\
\coun (K_i)=\coun (K_i^{-1})=1,\qquad \coun (E_i)=\coun (F_i)=0,\\
\antip (K_i)=K_i^{-1},\quad \antip (K_i^{-1})=K_i,\quad
\antip (E_i)=-E_iK_i^{-1},\quad \antip (F_i)=-K_iF_i
\end{gather}
for $i=1,2,\ldots ,n$.

For $i=1,2,\ldots ,n$ and $j=1,2,\ldots ,N$ set
\begin{align}
\lambda _{i,j}&:=
\begin{cases}
q^{-2\delta ^j_i+2\delta ^j_{i+1}-2\delta ^j_{i'-1}+2\delta ^j_{i'}},&
i=1,2,\ldots ,n-1,\\
q^{-2\delta ^j_n+2\delta ^j_{n+2}},& N=2n+1, i=n,\\
q^{-2\delta ^j_{n-1}-2\delta ^j_n+2\delta ^j_{n+1}+2\delta ^j_{n+2}},&
N=2n, i=n.
\end{cases}
\end{align}
Observe that $\lambda _{i,j}\lambda _{i,j'}=1$ for all $i,j$ and
\begin{align}\label{eq-alambda}
q_i^{a_{ij}}&=\begin{cases}
\lambda _{i,j+1}\lambda _{i,j'},& i\leq n, j<n,\\
\lambda _{i,n+1}\lambda _{i,n+2},& i\leq n, j=n.
\end{cases}
\end{align}
{}From \cite[Section 8.4.1]{b-KS}
we recall the formulas for the vector representation $T_1$ of $\Uqso N$
with highest weight $(1+\delta _{N3},0,\ldots ,0)$ and highest weight vector
$\mathbf{e}_N=\gcl _N$ with respect to the simple roots $\alpha _1,\alpha _2,
\ldots ,\alpha _n$.
\begin{equation}\label{eq-fundrep}
\begin{gathered}
T_1(E_iK_i^{-1})=q^2E_{i+1,i}-q^2E_{i',i'-1},\quad
T_1(F_i)=E_{i,i+1}-E_{i'-1,i'},\\
T_1(K_i)=D_i^{-1}D_{i+1}D^{-1}_{i'-1}D_{i'}\quad \text{for $N=2n+\epsilon $,
$i=1,2,\ldots ,n-1$,}\\
T_1(E_nK_n^{-1})=q'(q^2E_{n+1,n}-qE_{n+2,n+1}),\quad
T_1(F_n)=q'(q^2E_{n,n+1}-q^{-1}E_{n+1,n+2}),\\
T_1(K_n)=D_n^{-1}D_{n+2}\quad \text{for $N=2n+1$,}\\
T_1(E_nK_n^{-1})=q^2E_{n+1,n-1}-q^2E_{n+2,n},\quad
T_1(F_n)=E_{n-1,n+1}-E_{n,n+2},\\
T_1(K_n)=D_{n-1}^{-1}D_n^{-1}D_{n+1}D_{n+2}\quad \text{for $N=2n$.}
\end{gathered}
\end{equation}
Here we used the notation $q'=\qp ^{1/2}$ and
$D_j=\sum _{k=1}^Nq^{2\delta ^k_j}E_{k,k}$, $j=1,2,\ldots ,n$,
and $E_{k,l}$, $k,l=1,2,\ldots ,N$, are the matrix units.
The aim of this section is to prove the following theorems.

\begin{thm}\label{t-homuqgcl}
There exists an algebra map $\pi :\Uqso N\to \Cl $, $N=2n+\epsilon $,
such that
{\allowdisplaybreaks
\begin{align}
\pi (K_i)&=\sum _{l=0}^{i+1}
\sum _\dind{1\leq j_1<j_2<}{\cdots <j_l\leq i+1}\left( \prod _{r=1}^l
\frac{\lambda _{i,j_r}-q^{4l-4r}}{c^2\qp q^{N+1-2j_r}}\right)
\gcl _{j_1}\gcl _{j_2}\cdots \gcl _{j_l}\gcl _{j'_l}\cdots \gcl _{j'_1},\\
\pi (K_i^{-1})&=\sum _{l=0}^{i+1}
\sum _\dind{1\leq j_1<j_2<}{\cdots <j_l\leq i+1}\left( \prod _{r=1}^l
\frac{\lambda ^{-1}_{i,j_r}-q^{4l-4r}}{c^2\qp q^{N+1-2j_r}}\right)
\gcl _{j_1}\gcl _{j_2}\cdots \gcl _{j_l}\gcl _{j'_l}\cdots \gcl _{j'_1},\\
\pi (E_iK_i^{-1})&=\frac{q^{2i+1-N}}{c^2\qp }\sum _{l=0}^{i-1}
\sum _\dind{1\leq j_1<j_2<}{\cdots <j_l\leq i-1}\left( \prod _{r=1}^l
\frac{1-q^{4r}}{c^2\qp q^{N+1-2j_r}}\right)
\gcl _{j_1}\cdots \gcl _{j_l}\gcl _{i+1}\gcl _{i'}\gcl _{j'_l}
\cdots \gcl _{j'_1},\\
\pi (F_i)&=\frac{q^{2i+1-N}}{c^2\qp }\sum _{l=0}^{i-1}
\sum _\dind{1\leq j_1<j_2<}{\cdots <j_l\leq i-1}\left( \prod _{r=1}^l
\frac{1-q^{4r}}{c^2\qp q^{N+1-2j_r}}\right)
\gcl _{j_1}\cdots \gcl _{j_l}\gcl _i\gcl _{i'-1}\gcl _{j'_l}
\cdots \gcl _{j'_1}\\
\intertext{for $i=1,2,\ldots ,n-1$,}
\pi (K_n)&=q^{2-\epsilon }\sum _{l=0}^n
\sum _\dind{1\leq j_1<j_2<}{\cdots <j_l\leq n}\left( \prod _{r=1}^l
\frac{\lambda _{n,j_r}-q^{4l-4r}}{c^2\qp q^{N+1-2j_r}}\right)
\gcl _{j_1}\gcl _{j_2}\cdots \gcl _{j_l}\gcl _{j'_l}\cdots \gcl _{j'_1},\\
\pi (K_n^{-1})&=q^{-2+\epsilon }\sum _{l=0}^n
\sum _\dind{1\leq j_1<j_2<}{\cdots <j_l\leq n}\left( \prod _{r=1}^l
\frac{\lambda ^{-1}_{n,j_r}-q^{4l-4r}}{c^2\qp q^{N+1-2j_r}}\right)
\gcl _{j_1}\gcl _{j_2}\cdots \gcl _{j_l}\gcl _{j'_l}\cdots \gcl _{j'_1}\\
\intertext{for $N=2n+\epsilon $,}
\pi (E_nK_n^{-1})&=\frac{q'}{c^2\qp }\sum _{l=0}^{n-1}
\sum _\dind{1\leq j_1<j_2<}{\cdots <j_l\leq n-1}\left( \prod _{r=1}^l
\frac{1-q^{4r}}{c^2\qp q^{N+1-2j_r}}\right)
\gcl _{j_1}\cdots \gcl _{j_l}\gcl _{n+1}\gcl _{n+2}\gcl _{j'_l}
\cdots \gcl _{j'_1},\\
\pi (F_n)&=\frac{q^{-1}q'}{c^2\qp }\sum _{l=0}^{n-1}
\sum _\dind{1\leq j_1<j_2<}{\cdots <j_l\leq n-1}\left( \prod _{r=1}^l
\frac{1-q^{4r}}{c^2\qp q^{N+1-2j_r}}\right)
\gcl _{j_1}\cdots \gcl _{j_l}\gcl _n\gcl _{n+1}\gcl _{j'_l}
\cdots \gcl _{j'_1}\\
\intertext{for $N=2n+1$ and}
\pi (E_nK_n^{-1})&=\frac{q^{-1}}{c^2\qp }\sum _{l=0}^{n-2}
\sum _\dind{1\leq j_1<j_2<}{\cdots <j_l\leq n-2}\left( \prod _{r=1}^l
\frac{1-q^{4r}}{c^2\qp q^{N+1-2j_r}}\right)
\gcl _{j_1}\cdots \gcl _{j_l}\gcl _{n+1}\gcl _{n+2}\gcl _{j'_l}
\cdots \gcl _{j'_1},\\
\pi (F_n)&=\frac{q^{-1}}{c^2\qp }\sum _{l=0}^{n-2}
\sum _\dind{1\leq j_1<j_2<}{\cdots <j_l\leq n-2}\left( \prod _{r=1}^l
\frac{1-q^{4r}}{c^2\qp q^{N+1-2j_r}}\right)
\gcl _{j_1}\cdots \gcl _{j_l}\gcl _{n-1}\gcl _n\gcl _{j'_l}
\cdots \gcl _{j'_1},
\end{align}
}
for $N=2n$.
\end{thm}

\begin{thm}\label{t-uqgreps}
(i) The vector space $V =\comp \{\gcl _i\,|\,i=1,2,\ldots ,N\}\subset
\Cl $ is invariant under the left adjoint action
$\pi _\mathrm{ad}:\Uqso N\to \End{\Cl }$ of $\Uqso N$, where
\begin{align}\label{eq-piad}
(\pi _\mathrm{ad}f)(v)=\pi (f_{(1)})\,v\,\pi (\antip (f_{(2)})),\qquad
f\in \Uqso N,\quad  v\in \Cl.
\end{align}
The representation $\pi _\mathrm{ad}: \Uqso N\to \End{V}$ of $\Uqso N$ is
isomorphic to the vector representation
with highest weight $(1+\delta _{N3},0,\ldots ,0)$ and highest weight vector
$\gcl _N$.\\
(ii) For $N=2n+1$ both mappings
$\pi ^\nu :\Uqso{2n+1}\to \End{\mc{I}^\nu _\mathrm{L}}$, $\nu =-1,1$,
$(\pi ^\nu f)(\psi ):=\pi (f)\psi $ for $f\in \Uqso{2n+1}$,
$\psi \in \mc{I}^\nu _\mathrm{L}$, are isomorphic to the
spin representation of $\Uqso{2n+1}$ with highest weight $(0,0,\ldots ,0,1)$
and highest weight vector $\varphi ^\nu _{2n+1}$.\\
(iii) For $N=2n$ the mapping
$\pi ^\nu :\Uqso{2n}\to \End{\Cls[2n]{\nu }\varphi ^1_{2n}}$, where
$\nu =+,-$, $(\pi ^\nu f)(\psi ):=\pi (f)\psi $ for $f\in \Uqso{2n}$,
$\psi \in \Cls[2n]{\nu }\varphi ^1_{2n}$, is isomorphic to the
spin representation of $\Uqso{2n}$ with highest weight $(0,0,\ldots ,0,0,1)$
and $(0,0,\ldots ,0,1,0)$ and highest weight vector
$\varphi ^1_{2n}$ and $\gcl _n\varphi ^1_{2n}$ for $\nu=+$ and $\nu =-$,
respectively.
\end{thm}

\begin{bem}
In Proposition 3.1.3 and Theorem 3.1.2 in \cite{inp-DingFrenkel94}
the existence of a homomorphism from the algebra $U_q(\mathfrak{L}^\pm )$
(which is closely related to $\Uqso N$, see \cite[Sect.~8.5]{b-KS})
to $\Cl $ was proved. Theorem \ref{t-homuqgcl} gives explicit
formulas in terms of the generators of $\Uqso N$.
\end{bem}

In order to prove these assertions we have to study the structure
of $\Cl $ in more detail.

\begin{lemma}\label{l-munot0}
Let $\mu _i\in \comp $, $i=1,2,\ldots ,N$, and let $\mu _j=0$ for
some $j\in \{1,2,\ldots ,N\}$.Then there exists no
$z\in \Cl $, $z\not=0$, such that $z\gcl _i=\mu _i\gcl _iz$
for all $i$.
\end{lemma}

\begin{bew}
Suppose that $z\in \Cl $, $z\not= 0$, such that
$z\gcl _i=\mu _i\gcl _iz$ for $i=1,2,\ldots ,N$. This yields in turn
$z\Cl \subset \Cl z$. Further, $\mu _j=0$ gives $z\gcl _j=0$ and
$\Cl z\gcl _j=0$. Hence $\mc{J}\gcl _j=0$, where $\mc{J}$ denotes
the two-sided ideal of $\Cl $ generated by $z$.
Since $z\in \mc{J}$, we have $\mc{J}\not= \{0\}$.
Further, $\gcl _j\not= 0$ and $\mc{J}\gcl _j=0$ imply $1\notin \mc{J}$
and therefore $\mc{J}\not=\Cl $.
By Propositions \ref{s-evensimple} and \ref{s-oddsemis} we conclude
that $N=2n+1$ and $\mc{J}=\mc{I}_\eta $, $\eta \in \{+,-\}$.
The relation $\rho _\eta \in \mc{J}$ and equation (\ref{eq-rhogamma}) give
$\rho :=\gcl _1\gcl _2\cdots \gcl _n(\eta c+\gcl _{n+1})\in \mc{J}$.
Hence $\rho \gcl _j=0$ implies $j\leq n$. Applying the algebra
antiautomorphism $\tau $ to (\ref{eq-rhogamma}), from
$\tau (\rho _\eta )=\rho _\eta $, $\tau (\mc{I}_\eta )=\mc{I}_\eta $ we get
$\rho ':=(\eta c+\gcl _{n+1})\gcl _{n+2}\gcl _{n+3}\cdots \gcl _N\in \mc{J}$.
Now $\rho '\gcl _j=0$ is fulfilled only for $j\geq n+2$.
This is a contradiction.
\end{bew}

\begin{lemma}\label{l-zkhomog}
Let $k\in \mathbb{N}$, $k\leq n$, and let $\mu ,\mu '\in \compx $.
If $z\in \Cl $, $z\not=0$, satisfies
$z\gcl _k=\mu \gcl _kz$ and $z\gcl _{k'}=\mu '\gcl _{k'}z$,
then $z$ is homogeneous with respect to $\partial _k$ and we have
$\partial _k(z)=0$.
\end{lemma}

\begin{bew}
There exist unique elements $z^{(r)}\in \Cl $, $r=-1,0,1$, such that
$\partial _k(z^{(r)})=r$ and $z=z^{(-1)}+z^{(0)}+z^{(1)}$.
The action of the algebra automorphism
$\gcl _i\mapsto q^{\delta ^i_k-\delta ^i_{k'}}\gcl _i$, $i=1,2,\ldots ,N$,
and its square on the equation $z\gcl _k-\mu \gcl _kz=0$ give
\begin{align}
(q^{-1}z^{(-1)}+z^{(0)}+qz^{(1)})\gcl _k&=
\mu \gcl _k(q^{-1}z^{(-1)}+z^{(0)}+qz^{(1)}),\\
(q^{-2}z^{(-1)}+z^{(0)}+q^2z^{(1)})\gcl _k&=
\mu \gcl _k(q^{-2}z^{(-1)}+z^{(0)}+q^2z^{(1)}).
\end{align}
Together with $z\gcl _k=\mu \gcl _kz$ we obtain $z^{(r)}\gcl _k=
\mu \gcl _kz^{(r)}$ for $r=-1,0,1$. Similarly,
$z^{(r)}\gcl _{k'}=\mu '\gcl _{k'}z^{(r)}$ for $r=-1,0,1$.
We will show that $z^{(-1)}=z^{(1)}=0$.

By Proposition \ref{s-clbasis}, $z^{(-1)}$ can uniquely be written as
$z^{(-1)}=\sum _{I_N}\lambda _{I_N}\gcl ^{I_N}$.
Since $\partial _k(z^{(-1)})=-1$, for all multiindices $I_N$ with
$\lambda _{I_N}\not=0$ we have $i_k=0$, $i_{k'}=1$.
Let us fix a multiindex $J_N$ in the set of those $I_N$ for which
$\lambda _{I_N}\not=0$ and the sum $j_1+\cdots +j_{k-1}$ is minimal, say $m$.
Consider the element
\begin{align*}
\gcl _1^{1-j_1}\gcl _2^{1-j_2}\cdots \gcl _{k-1}^{1-j_{k-1}}
\gcl _k(z^{(-1)}\gcl _k-\mu \gcl _kz^{(-1)}).
\end{align*}
Since $\gcl _k^2=0$, the
second summand vanishes. Using the relations (\ref{eq-cliff}) the first
summand becomes
\begin{align*}
&\gcl _1^{1-j_1}\gcl _2^{1-j_2}\cdots \gcl _{k-1}^{1-j_{k-1}}\gcl _kz^{(-1)}
\gcl _k=\\
&\qquad =\gcl _1^{1-j_1}\gcl _2^{1-j_2}\cdots \gcl _{k-1}^{1-j_{k-1}}\gcl _k
\sum _{I_N}\lambda _{I_N}\gcl _1^{i_1}\gcl _2^{i_2}\cdots 
\gcl _N^{i_N}\gcl _k=\\
&\qquad =\gcl _1^{1-j_1}\gcl _2^{1-j_2}\cdots \gcl _{k-1}^{1-j_{k-1}}\gcl _k
\sum _{I_N}\delta ^{i_k}_0\delta ^{i_{k'}}_1(-q^2)^{i_{k'+1}+\cdots +i_N}
\lambda _{I_N}\gcl _1^{i_1}\gcl _2^{i_2}\cdots \gcl _{k'-1}^{i_{k'-1}}\times \\
&\qquad \phantom{=}\times
\left(-\gcl _k\gcl _{k'}+\qm \qp \sum _{l=1}^{k-1}q^{2l-2k+2}\gcl _l\gcl _{l'}
+c^2q^{N-2k+1}\qp \right)\gcl _{k'+1}^{i_{k'+1}}\cdots \gcl _N^{i_N}.
\end{align*}
Since $\gcl _k^2=0$, the first summand inside the large brackets vanishes
for all $I_N$.
Recall that if $i_k=0$ and $i_{k'}=1$ then $i_1+i_2+\ldots +i_{k-1}\geq m$.
Hence for the number of factors $\gcl _t$, $1\leq t<k$,
in the second expression inside of the brackets we obtain
\begin{align*}
(1-j_1)&+(1-j_2)+\ldots +(1-j_{k-1})+i_1+i_2+\cdots +i_{k-1}+1\geq \\
&\geq k-1-j_1-j_2-\cdots -j_{k-1}+m+1=k-1-m+m+1=k.
\end{align*}
This means that these summands are zero as well. Therefore we are left with
\begin{align*}
\gcl _1^{1-j_1}\gcl _2^{1-j_2}&\cdots \gcl _{k-1}^{1-j_{k-1}}\gcl _kz^{(-1)}
\gcl _k=c^2q^{N-2k+1}\qp \sum _{I_N}\delta ^{i_k}_0\delta ^{i_{k'}}_1
(-q^2)^{n(I_N,J_N)}\lambda _{I_N}\times \\
&\times \gcl _1^{1-j_1+i_1}\gcl _2^{1-j_2+i_2}\cdots 
\gcl _{k-1}^{1-j_{k-1}+i_{k-1}}
\gcl _k \gcl _{k+1}^{i_{k+1}}\cdots \gcl _{k'-1}^{i_{k'-1}}
\gcl _{k'+1}^{i_{k'+1}}\cdots \gcl _N^{i_N}
\end{align*}
with well-defined integer numbers $n(I_N,J_N)$. The appearing nonzero summands
are linearly independent because of Proposition \ref{s-clbasis}. Then
$z^{(-1)}\gcl _k-\mu \gcl _kz^{(-1)}=0$ implies that the above sum must be
zero. In particular the summand for $I_N=J_N$ is zero if and only if
$\lambda _{J_N}=0$ which is a contradiction. Hence $z^{(-1)}=0$.

Observe that $\partial _k(\tau (z^{(r)}))=-r$ and
\begin{align}\label{eq-tauzgcl}
\tau (z)\gcl _k=(\mu ')^{-1}\gcl _k\tau (z).
\end{align}
Replacing $\mu $ by $(\mu ')^{-1}$ and $z^{(-1)}$ by $\tau (z^{(1)})$ in the
above proof, we obtain $\tau (z^{(1)})= 0$ and hence $z^{(1)}=0$.
Therefore $z=z^{(0)}$.
\end{bew}

\begin{lemma}\label{l-muvector}
Let $\mu _i\in \compx $, $i=1,2,\ldots ,N$, and let $z\not =0$ be an element
in $\Cl $ such that $z\gcl _i=\mu _i\gcl _iz$ for all
$i=1,2,\ldots ,N$. Then $\mu _i\mu _{i'}=1$ and $\mu _j=1$ for $j=j'$.
\end{lemma}

\begin{bew}
{}From Lemma \ref{l-zkhomog} we conclude that $\partial _k(z)=0$ for all
$k=1,2,\ldots ,n$. Hence $\tau (z)=z$. The equations
$\mu _k\gcl _kz-z\gcl _k=0$ and (\ref{eq-tauzgcl}) give
\begin{align}
0&=\mu _k\gcl _kz-\mu _k\gcl _k\tau (z)=z\gcl _k-\mu _k\mu _{k'}\tau (z)\gcl _k
=(1-\mu _k\mu _{k'})z\gcl _k.
\end{align}
Since $z\gcl _k\not=0$ by Lemma \ref{l-munot0}, we get $\mu _k\mu _{k'}=1$.
Finally, if $N=2n+1$ then the elements $z$ and $\gcl _{n+1}$ commute by
Proposition \ref{s-clzcomm}.
\end{bew}

\begin{lemma}\label{l-partZcrit}
Let $\mu =(\mu _1,\mu _2,\ldots ,\mu _k)\in (\compx )^k$, $1\leq k\leq n$.
Then the following assertions are equivalent:\\
(i) $z\gcl _i=\mu _i\gcl _iz$, $\tau (z)\gcl _i=\mu _i\gcl _i\tau (z)$
for $i=1,2,\ldots ,k$,\\
(ii) $\partial _i(z)=0$ and
\begin{align*}
z\gcl _1\gcl _2\cdots \gcl _i&=\mu _1\mu _2\cdots \mu _i
\gcl _1\gcl _2\cdots \gcl _iz,& 
\tau (z)\gcl _1\gcl _2\cdots \gcl _i&=\mu _1\mu _2\cdots \mu _i
\gcl _1\gcl _2\cdots \gcl _i\tau (z)
\end{align*}
for $i=1,2,\ldots ,k$,\\
(iii) $\partial _i(z)=0$ and
$\gcl _1\gcl _2\cdots \gcl _{i-1}(z\gcl _i-\mu _i\gcl _iz)=
\gcl _1\gcl _2\cdots \gcl _{i-1}(\tau (z)\gcl _i
-\mu _i\gcl _i\tau (z))=0$ for $i=1,2,\ldots ,k$.
\end{lemma}

\begin{bew}
Applying the algebra antiautomorphism $\tau $ of $\Cl $ to the second equation
of (i) and using $\tau ^2=\id $
the conclusion (i)$\to $(iii) follows immediately from Lemma \ref{l-zkhomog}.
For the proof of (iii)$\to $(ii) we take the sequence of equations
\begin{align*}
\mu _1\mu _2&\cdots \mu _j\gcl _1\gcl _2\cdots \gcl _jz'=
\mu _1\mu _2\cdots \mu _{j-1}(\gcl _1\gcl _2\cdots \gcl _{j-1}z'\gcl _j)=\\
&=\mu _1\cdots \mu _{j-2}(\gcl _1\cdots \gcl _{j-2}z'\gcl _{j-1})
\gcl _j=\cdots =\mu _1\gcl _1z'\gcl _2\cdots \gcl _j
=z'\gcl _1\gcl _2\cdots \gcl _j,
\end{align*}
where the equations of (iii) for $i=j,j-1,\ldots ,2,1$ and $z'=z,\tau(z)$
are used.

We have to prove (ii)$\to $(i).
Suppose that $z\in \Cl $, $z\not=0$, fulfills (ii). We prove by induction that
$z\gcl _i=\mu _i\gcl _iz$ and $\tau (z)\gcl _i=\mu _i\gcl _i\tau (z)$
for $i=1,2,\ldots ,k$.
Our induction assumptions are the equations in (ii) for fixed
$i:=m<k$ and the equations $z\gcl _j=\mu _j\gcl _jz$,
$z\gcl _{j'}=\mu _j^{-1}\gcl _{j'}z$ for all $j=1,2,\ldots ,m-1$.
Observe that the last requirements are empty at the beginning $m=1$ of the
induction. We have to prove $z\gcl _m=\mu _m\gcl _mz$ and
$\tau (z)\gcl _m=\mu _m\gcl _m\tau (z)$. Applying $\tau $ to the last equation
we then obtain $z\gcl _{m'}=\mu _m^{-1}\gcl _{m'}z$ which is needed for the
next induction step.

Now we prove not only $z\gcl _m=\mu _m\gcl _mz$ but all equations
\begin{align}\label{eq-jinduction}
z\gcl _1\gcl _2\cdots \gcl _j\gcl _m=\mu _1\mu _2\cdots \mu _j\mu _m
\gcl _1\gcl _2\cdots \gcl _j\gcl _mz,
\end{align}
$j=0,1,2,\ldots ,m-1$ by inverse induction on $j$. For $j=m-1$ this is the
induction hypothesis and for $j=0$ the equation we want to prove.
For the induction step we use the formula $z\gcl _{j'}=\mu _j^{-1}\gcl _{j'}z$
which is known since $j<m$. By (\ref{eq-cliff}) we compute
\begin{align*}
&z\gcl _1\gcl _2\cdots \gcl _{j-1}\gcl _m=
z\left(\frac{-q^{2j+1-N}}{\qp c^2}\gcl _1\gcl _2\cdots \gcl _j\gcl _m
\gcl _{j'}-(-1)^j\frac{q^{1-N}}{\qp c^2}\gcl _{j'}\gcl _1\gcl _2\cdots
\gcl _j\gcl _m\right)\\
&=\mu _1\mu _2\cdots \mu _j\mu _m\mu _j^{-1}
\left(\frac{-q^{2j+1-N}}{\qp c^2}\gcl _1\gcl _2\cdots \gcl _j
\gcl _m\gcl _{j'}-(-1)^j\frac{q^{1-N}}{\qp c^2}\gcl _{j'}\gcl _1\gcl _2
\cdots \gcl _j\gcl _m\right)z\\
&=\mu _1\mu _2\cdots \mu _{j-1}\mu _m\gcl _1\gcl _2\cdots \gcl _{j-1}\gcl _mz.
\end{align*}
Changing the role of $z$ and $\tau (z)$ the remaining assertion
$\tau (z)\gcl _m=\mu _m\gcl _m\tau (z)$ can be shown in the same way.
\end{bew}

Taking Lemma \ref{l-munot0} and Lemma \ref{l-muvector} into account
we define complex vector spaces $\mc{Z}^\mu _N$ for
$\mu =(\mu _1,\mu _2,\ldots ,\mu _n)\in (\compx )^n$. Set
$\mu _{i'}=\mu _i^{-1}$ for $i\in \{1,2,\ldots ,N\}$, $i>i'$, and
$\mu _i=1$ for $i=i'$. Then we define
\begin{align}
\mc{Z}^\mu _N=\{z\in \Cl \,|\,
z\gcl _i=\mu _i\gcl _iz,\quad i=1,2,\ldots ,N\}.
\end{align}
Similarly to the proof of Lemma \ref{l-zkhomog}
one can show that $\mc{Z}^\mu _N$ splits into the direct sum
$\mc{Z}^{\mu ,0}_N\oplus \mc{Z}^{\mu ,1}_N$
of $\partial _0$-homogeneous components of even and odd degree, respectively.
Since there are no nonzero elements $z$ of $\Clo[2n]$ such that
$\partial _k(z)=0$ for all $k=1,2,\ldots ,n$, we have
$\mc{Z}^{\mu ,1}_{2n}=\{0\}$.

\begin{lemma}\label{l-fullZcrit}
Let $\mu =(\mu _1,\mu _2,\ldots ,\mu _n)\in (\compx )^n$.
Then $z\in \mc{Z}^\mu _N$ if and only if
$\partial _k(z)=0$ and
$\gcl _1\gcl _2\cdots \gcl _{k-1}(z\gcl _k-\mu _k\gcl _kz)=0$
for all $k=1,2,\ldots ,n$.
\end{lemma}

\begin{bew}
The only if part of the lemma follows from Lemma
\ref{l-zkhomog} and the definition of $\mc{Z}^\mu _N$.
On the other hand, suppose that 
$\partial _k(z)=0$ and
$\gcl _1\gcl _2\cdots \gcl _{k-1}(z\gcl _k-\mu _k\gcl _kz)=0$
for all $k=1,2,\ldots ,n$.
Then $\tau (z)=z$ and hence we can apply
Lemma~\ref{l-partZcrit}(iii)$\to $(i) with $k=n$.
This means $z\gcl _i=\mu _i\gcl _iz$ and $z\gcl _{i'}=\mu _i^{-1}\gcl _{i'}z$.
The remaining assertion $z\gcl _{n+1}=\gcl _{n+1}z$ follows from
$\{z,\gcl _{n+1}\}\subset \Clz[2n+1]$ and Proposition \ref{s-clzcomm}.
\end{bew}

\begin{satz}\label{s-mucenter}
Let $\mu =(\mu _1,\mu _2,\ldots ,\mu _n)\in (\compx )^n$.
The vector spaces $\mc{Z}^{\mu ,0}_{2n},\mc{Z}^{\mu ,0}_{2n+1}$
and $\mc{Z}^{\mu ,1}_{2n+1}$ are one-dimensional
and generated by the element
\begin{align}
\notag
z_{\mu ,\hat{\eta }}=\sum _{k=0}^n
\sum _\dind{1\leq i_1<i_2<}{\cdots <i_k\leq n}\left(
\prod _{r=1}^k\frac{\mu _{i_r}-(-q^2)^{\hat{\eta }}q^{4k-4r}}{c^2\qp
q^{N+1-2i_r}}\right)
\gcl _{i_1}\gcl _{i_2}\cdots \gcl _{i_k}\gcl _{n+1}^{\hat{\eta }}
\gcl _{i'_k}\cdots \gcl _{i'_2}\gcl _{i'_1},
\end{align}
where $\hat{\eta }=0$ for $\mc{Z}^{\mu ,0}_{2n}$ and $\mc{Z}^{\mu ,0}_{2n+1}$,
and $\hat{\eta }=1$ for $\mc{Z}^{\mu ,1}_{2n+1}$.
\end{satz}

\begin{folg}\label{f-z0z1}
The elements $z_0\in \Cl[2n]$ and
$z_1\in \Clo[2n+1]$, given by
\begin{align}
z_0&=\sum _{k=0}^n(-1)^k\prod _{l=1}^k
\frac{q^{2l-2}+q^{-2l+2}}{q\qp c^2}
\sum _\dind{J_n}{\ell (J_n)=k} q^{2L(J_n)-k(2n-k+1)}\gcl _1^{j_1}\gcl _2^{j_2}
\cdots \gcl _n^{j_n}\gcl _{n+1}^{j_n}\cdots \gcl _N^{j_1},\\
z_1&=\frac{1}{c}\sum _{k=0}^n\prod _{l=1}^k
\frac{q^{2l-1}+q^{-2l+1}}{q\qp c^2}
\sum _\dind{J_n}{\ell (J_n)=k} q^{2L(J_n)-k(2n-k+1)}\gcl _1^{j_1}\gcl _2^{j_2}
\cdots \gcl _n^{j_n}\gcl _{n+1}\gcl _{n'}^{j_n}\cdots \gcl _N^{j_1},
\end{align}
satisfy the equations $z_0\gcl _i=-\gcl _iz_0$ and
$z_1\gcl _i=\gcl _iz_1$, respectively,
for all $i=1,2,\ldots ,N$. Moreover, they are unique with this property
up to a scalar factor.
\end{folg}

\begin{bew}[ of the proposition]
Let $z\in \mc{Z}^{\mu ,\hat{\eta }}_N$, $z\not= 0$. We use
Lemma \ref{l-fullZcrit} to determine the explicit form of $z$.
The requirement $\partial _k(z)=0$ for all
$k=1,2,\ldots ,n$ is equivalent to the fact that $z$ takes the form
\begin{align*}
z&=\sum _{I_n}\lambda _{I_n}\gcl _1^{i_1}\gcl _2^{i_2}\cdots \gcl _n^{i_n}
\gcl _{n+1}^{\hat{\eta }}\gcl _{n'}^{i_n}\cdots \gcl _N^{i_1}.
\end{align*}
If $I_n$ is a multiindex with $i_k=0$ for some $k\in \{1,2,\ldots ,n\}$,
then let $I_n^k$ denote the multiindex $(j_1,j_2,\ldots ,j_n)$ such that
$j_l=i_l+\delta _{lk}$, $l=1,2,\ldots ,n$.
Let us compute the left hand side of the remaining
equations in Lemma \ref{l-fullZcrit}. By (\ref{eq-cliff}) we obtain
\begin{align}\label{eq-gclz}
\gcl _1\gcl _2&\cdots \gcl _{k-1}(z\gcl _k-\mu _k\gcl _kz)=\\
\notag
&=\sum _{I_n,i_l=0 \text{ for }l<k}\lambda _{I_n}
\gcl _1\gcl _2\cdots \gcl _{k-1}(\gcl _k^{i_k}\gcl _{k+1}^{i_{k+1}}
\cdots \gcl _n^{i_n}\gcl _{n+1}^{\hat{\eta }}\gcl _{n'}^{i_n}\cdots
\gcl _{k'}^{i_k}\gcl _k\\
\notag &
\phantom{=\sum _{I_n,i_l=0 \text{ for }l<k}\lambda _{I_n}
\gcl _1\gcl _2\cdots \gcl _{k-1}(}
-\mu _k\gcl _k\gcl _k^{i_k}\gcl _{k+1}^{i_{k+1}}\cdots \gcl _n^{i_n}
\gcl _{n+1}^{\hat{\eta }}\gcl _{n'}^{i_n}\cdots \gcl _{k'}^{i_k})\\
\notag
&=\sum _{I_n,i_l=0 \text{ for }l\leq k}\gcl _1\gcl _2\cdots \gcl _{k-1}
(\lambda _{I^k_n}\gcl _k\gcl _{k+1}^{i_{k+1}}\cdots \gcl _{k'-1}^{i_{k+1}}
\gcl _{k'}\gcl _k+\lambda _{I_n}\gcl _{k+1}^{i_{k+1}}\cdots
\gcl _{k'-1}^{i_{k+1}}\gcl _k\\
\notag
&\phantom{=\sum _{I_n,i_l=0 \text{ for }l\leq k}\gcl _1\gcl _2\cdots
\gcl _{k-1}(}
-\mu _k\lambda _{I_n}\gcl _k\gcl _{k+1}^{i_{k+1}}\cdots
\gcl _{k'-1}^{i_{k+1}})\\
\notag
&=\sum _{I_n,i_l=0 \text{ for }l\leq k}\gcl _1\gcl _2\cdots \gcl _{k-1}
(\lambda _{I^k_n}c^2q^{N-2k+1}\qp +\lambda _{I_n}(-q^2)^{2i_{k+1}+2i_{k+2}+
\cdots +2i_n+{\hat{\eta }}}\\
\notag
&\phantom{=\sum _{I_n,i_l=0 \text{ for }l\leq k}\gcl _1\gcl _2\cdots
\gcl _{k-1}(}
-\mu _k\lambda _{I_n})\gcl _k\gcl _{k+1}^{i_{k+1}}\cdots
\gcl _{k'-1}^{i_{k+1}}.
\end{align}
By Proposition \ref{s-clbasis} the summands with nonzero coefficient in the
last expression are linearly independent. Hence 
$\gcl _1\gcl _2\cdots \gcl _{k-1}(z\gcl _k-\mu _k\gcl _kz)=0$ if and only if
\begin{align}\label{eq-lambdarec}
\lambda _{I^k_n}&=\frac{\mu _k-(-q^2)^{2i_{k+1}+2i_{k+2}+\cdots 
+2i_n+{\hat{\eta }}}}{c^2q^{N-2k+1}\qp }\lambda _{I_n}
\end{align}
for all $k=1,2,\ldots ,n$. This yields $z\in \comp z_{\mu ,\hat{\eta }}$.
\end{bew}

Recall the definition of $z_0$ and $z_1$ in Corollary \ref{f-z0z1}.

\begin{folg}
We have $z_0^2=1$ in $\Cl[2n]$ and $z_1^2=1$ in
$\Cl[2n+1]$.
\end{folg}

\begin{bew}
By definition $z_\epsilon $ is a $\partial _0$-homogeneous
element in $\Cl[2n+\epsilon ]$.
Hence $z_\epsilon ^2\in \Cle[2n+\epsilon ]$ and by Corollary \ref{f-z0z1}
it belongs to the center
of $\Cl[2n+\epsilon ]$. By Proposition \ref{s-mucenter} with
$\mu :=(1,1,\ldots ,1)$, this implies
$z_\epsilon ^2=\lambda 1$ for some $\lambda \in \comp $.
Multiplying this equation by $\varphi ^1_N$ from the right and taking
into account the formula
$z_\epsilon \varphi ^1_N=\varphi ^1_N$ we conclude that
$\lambda =1$.
\end{bew}

\begin{folg}
The element $z_1$ acts (by multiplication) with the factor $\eta $
on the left ideal
$\mc{I}^\eta _\mathrm{L}$ of $\Cl[2n+1]$. The element
$z_0$ acts with the factor $+1$ and
$-1$ on the subspaces $\Cle[2n]\varphi ^1_{2n}$
and $\Clo[2n]\varphi ^1_{2n}$
of the left ideal $\mc{I}^1_\mathrm{L}$ of $\Cl[2n]$,
respectively. 
\end{folg}

Now we turn to the proofs of Theorem \ref{t-homuqgcl} and
Theorem \ref{t-uqgreps}. Therein the following lemma plays a crucial role.

\begin{lemma}\label{l-pirels}
For $f=K_i,K_i^{-1},E_iK_i^{-1},F_i$, $i=1,2,\ldots ,n$, let $\tilde{f}$
denote the element of $\Cl $ given by the expression
at the right of $\pi (f)$ in Theorem \ref{t-homuqgcl}.
Then the elements $\tilde{f}$ satisfy the following equations in $\Cl $:
{\allowdisplaybreaks
\begin{gather}
\label{eq-Kgcl}
\tilde{K}_i\gcl _j=\lambda _{i,j}\gcl _j\tilde{K}_i,\qquad
\widetilde{K_i^{-1}}\gcl _j=\lambda ^{-1}_{i,j}\gcl _j\widetilde{K^{-1}_i},\\
\label{eq-Egcl}
\widetilde{E_kK_k^{-1}}\gcl _j=\begin{cases}
\gcl _j\widetilde{E_kK_k^{-1}} & \text{for $j\notin \{k,k+1,k'-1,k'\}$,}\\
q^{-2}\gcl _j\widetilde{E_kK_k^{-1}} & \text{for $j=k+1,k'$,}\\
q^2\gcl _j\widetilde{E_kK_k^{-1}}+q^2\gcl _{j+1} & \text{for $j=k$,}\\
q^2\gcl _j\widetilde{E_kK_k^{-1}}-q^2\gcl _{j+1} & \text{for $j=k'-1$,}
\end{cases}\\
\label{eq-Fgcl}
\tilde{F}_k\gcl _j=\begin{cases}
\gcl _j\tilde{F_k} & \text{for $j\notin \{k,k+1,k'-1,k'\}$,}\\
q^2\gcl _j\tilde{F_k} & \text{for $j=k,k'-1$,}\\
q^{-2}\gcl _j\tilde{F_k}+\gcl _{j-1} & \text{for $j=k+1$,}\\
q^{-2}\gcl _j\tilde{F_k}-\gcl _{j-1} & \text{for $j=k'$,}
\end{cases}\\
\label{eq-Engcl}
\widetilde{E_nK_n^{-1}}\gcl _j=\begin{cases}
\gcl _j\widetilde{E_nK_n^{-1}} & \text{for $j\notin \{n-1,n,n+1,n+2\}$,}\\
\gcl _j\widetilde{E_nK_n^{-1}} & \text{for $j=n-1$, $N=2n+1$,}\\
q^2\gcl _j\widetilde{E_nK_n^{-1}}+q^2q'\gcl _{j+1} & \text{for $j=n$,
 $N=2n+1$,}\\
\gcl _j\widetilde{E_nK_n^{-1}}-qq'\gcl _{j+1} & \text{for $j=n+1$, $N=2n+1$,}\\
q^{-2}\gcl _j\widetilde{E_nK_n^{-1}} & \text{for $j=n+2$, $N=2n+1$,}\\
q^2\gcl _j\widetilde{E_nK_n^{-1}}+q^2\gcl _{j+2} & \text{for $j=n-1$,
 $N=2n$,}\\
q^2\gcl _j\widetilde{E_nK_n^{-1}}-q^2\gcl _{j+2} & \text{for $j=n$, $N=2n$,}\\
q^{-2}\gcl _j\widetilde{E_nK_n^{-1}} & \text{for $j=n+1,n+2$, $N=2n$,}
\end{cases}\\
\label{eq-Fngcl}
\tilde{F_n}\gcl _j=\begin{cases}
\gcl _j\tilde{F_n} & \text{for $j\notin \{n-1,n,n+1,n+2\}$,}\\
\gcl _j\tilde{F_n} & \text{for $j=n-1$, $N=2n+1$,}\\
q^2\gcl _j\tilde{F_n} & \text{for $j=n$, $N=2n+1$,}\\
\gcl _j\tilde{F_n}+q'\gcl _{j-1} & \text{for $j=n+1$, $N=2n+1$,}\\
q^{-2}\gcl _j\tilde{F_n}-q^{-1}q'\gcl _{j-1} & \text{for $j=n+2$, $N=2n+1$,}\\
q^2\gcl _j\tilde{F_n} & \text{for $j=n-1,n$, $N=2n$,}\\
q^{-2}\gcl _j\tilde{F_n}+\gcl _{j-2} & \text{for $j=n+1$, $N=2n$,}\\
q^{-2}\gcl _j\tilde{F_n}-\gcl _{j-2} & \text{for $j=n+2$, $N=2n$,}
\end{cases}
\end{gather}
}
for all $i=1,2,\ldots ,n$ and $k=1,2,\ldots ,n-1$, where $q'=\qp ^{1/2}$.
\end{lemma}

\begin{bew}
Equation (\ref{eq-Kgcl}) follows in turn from Proposition \ref{s-mucenter}.
Observe that by definition $\widetilde{E_iK_i^{-1}}=
q^{\delta _{\epsilon 1}\delta _{in}}\tau (\tilde{F}_i)$ for all
$i=1,2,\ldots ,n$. Moreover, applying $\tau $ to the equations
(\ref{eq-Egcl}), (\ref{eq-Fgcl}), (\ref{eq-Engcl}) and (\ref{eq-Fngcl})
we obtain the formulas (\ref{eq-Fgcl}), (\ref{eq-Egcl}), (\ref{eq-Fngcl})
and (\ref{eq-Engcl}), respectively, for $j'$ instead of $j$.
Therefore it suffices to show (\ref{eq-Egcl})--(\ref{eq-Fngcl}) for $j\leq j'$.
Let us prove (\ref{eq-Egcl}) and (\ref{eq-Fgcl}). The proof of the equations
(\ref{eq-Engcl}) and (\ref{eq-Fngcl}) can be treated similarly.

If $j\in \{k+2,k+3,\ldots ,n+\epsilon \}$ then
$\tilde{F}_k\gcl _j=\gcl _j\tilde{F}_k$ and
$\widetilde{E_kK_k^{-1}}\gcl _j=\gcl _j\widetilde{E_kK_k^{-1}}$.
Indeed, in both $\tilde{F}_k$ and $\widetilde{E_kK_k^{-1}}$ there are
no summands containing factors $\gcl _l$, $k+1<l<k'-1$. Moreover, in each
summand there are
as many factors $\gcl _l$, $l\leq k+1$, as factors $\gcl _l$, $l\geq k'-1$.
{}From this and the first line of (\ref{eq-cliff}) it follows that each summand
of $\tilde{F}_k$ and of $\widetilde{E_kK_k^{-1}}$ commutes with $\gcl _j$.

If $j<k$ then we use Lemma \ref{l-partZcrit}(i)$\Leftrightarrow $(iii) with
$(\mu _1,\mu _2,\ldots ,\mu _j)=(1,1,\ldots ,1)$. Obviously,
$\partial _i(\widetilde{E_kK_k^{-1}})=\partial _i(\tilde{F}_k)=0$ for
$i\leq j$. Moreover,
\begin{align*}
&\gcl _1\gcl _2\cdots \gcl _{j-1}(\tilde{F}_k\gcl _j-\gcl _j\tilde{F}_k)=\\
&=\gcl _1\gcl _2\cdots \gcl _{j-1}
\sum _{I_n,i_l=0 \text{ for }l<j}\lambda _{I_n}(
\gcl _j^{i_j}\gcl _{j+1}^{i_{j+1}}\cdots \gcl _{k-1}^{i_{k-1}}\gcl _k
\gcl _{k'-1}\gcl _{k'+1}^{i_{k-1}}\cdots \gcl _{j'-1}^{i_{j+1}}\gcl _{j'}^{i_j}
\gcl _j\\
&\phantom{
=\gcl _1\gcl _2\cdots \gcl _{j-1}
\sum _{I_n,i_l=0 \text{ for }l<j}\lambda _{I_n}(
}
-\gcl _j\gcl _j^{i_j}\gcl _{j+1}^{i_{j+1}}\cdots \gcl _{k-1}^{i_{k-1}}\gcl _k
\gcl _{k'-1}\gcl _{k'+1}^{i_{k-1}}\cdots \gcl _{j'-1}^{i_{j+1}}\gcl _{j'}^{i_j}
)\\
&=\gcl _1\gcl _2\cdots \gcl _{j-1}\sum _{I_n,i_l=0 \text{ for }l\leq j}(
\lambda _{I^j_n}\gcl _j\gcl _{j+1}^{i_{j+1}}\cdots \gcl _{k-1}^{i_{k-1}}\gcl _k
\gcl _{k'-1}\gcl _{k'+1}^{i_{k-1}}\cdots \gcl _{j'-1}^{i_{j+1}}\gcl _{j'}
\gcl _j\\
&\phantom{
=\gcl _1\gcl _2\cdots \gcl _{j-1}\sum _{I_n,i_l=0 \text{ for }l<j}(
}
+\lambda _{I_n}\gcl _{j+1}^{i_{j+1}}\cdots \gcl _{k-1}^{i_{k-1}}\gcl _k
\gcl _{k'-1}\gcl _{k'+1}^{i_{k-1}}\cdots \gcl _{j'-1}^{i_{j+1}}\gcl _j
\\
&\phantom{
=\gcl _1\gcl _2\cdots \gcl _{j-1}\sum _{I_n,i_l=0 \text{ for }l\leq j}(
}
-\lambda _{I_n}\gcl _j\gcl _{j+1}^{i_{j+1}}\cdots \gcl _{k-1}^{i_{k-1}}\gcl _k
\gcl _{k'-1}\gcl _{k'+1}^{i_{k-1}}\cdots \gcl _{j'-1}^{i_{j+1}}
)\\
&=\sum _{I_n,i_l=0 \text{ for }l\leq j}\gcl _1\gcl _2\cdots \gcl _{j-1}(
\lambda _{I^j_n}c^2q^{N-2j+1}\qp +\lambda _{I_n}(-q^2)^{2i_{j+1}+\cdots
+2i_{k-1}+2}\\
&\phantom{
=\sum _{I_n,i_l=0 \text{ for }l\leq j}\gcl _1\gcl _2\cdots \gcl _{j-1}(
}
-\lambda _{I_n})\gcl _j\gcl _{j+1}^{i_{j+1}}\cdots \gcl _{k-1}^{i_{k-1}}\gcl _k
\gcl _{k'-1}\gcl _{k'+1}^{i_{k-1}}\cdots \gcl _{j'-1}^{i_{j+1}}.
\end{align*}
The last expression vanishes because of
$\displaystyle \lambda _{I_n^j}=\frac{1-q^{4\ell (I_n)+4}}{c^2q^{N+1-2j}\qp }$
for those coefficients $\lambda _{I_n}$ of $\tilde{F}_k$ for which
$i_l=0$ for $l\leq j$ and $l\geq k$. Hence
$\gcl _1\gcl _2\cdots \gcl _{j-1}(\tilde{F}_k\gcl _j-\gcl _j\tilde{F}_k)=0$.
Similarly one proves
$\gcl _1\gcl _2\cdots \gcl _{j-1}(\widetilde{E_kK_k^{-1}}\gcl _j-\gcl _j
\widetilde{E_kK_k^{-1}})=0$. Hence, by Lemma \ref{l-partZcrit}
we get (\ref{eq-Egcl}) and (\ref{eq-Fgcl}) for $j<k$.

In a completely similar way one checks that
\begin{equation}\label{eq-gclEFgcl}
\begin{gathered}
\gcl _1\gcl _2\cdots \gcl _{k-1}(\widetilde{E_kK_k^{-1}}\gcl _k-q^2
\gcl _k\widetilde{E_kK_k^{-1}}-q^2\gcl _{k+1})=0,\\
\gcl _1\gcl _2\cdots \gcl _{k-1}(\tilde{F}_k\gcl _{k+1}-q^{-2}\gcl _{k+1}
\tilde{F}_k-\gcl _k)=0.
\end{gathered}
\end{equation}
Multiplying these equations from the left and from the right by
$\gcl _{j'}$, $j=k-1,k-2,\ldots ,1$, respectively, using (\ref{eq-Egcl})
and (\ref{eq-Fgcl}) for $j>k'$, and taking the appropriate
linear combination of the resulting expressions, one can successively
cancel the leading factors $\gcl _j$, $j=k-1,k-2,\ldots ,1$ in
(\ref{eq-gclEFgcl}).
This proves (\ref{eq-Egcl}) for $j=k$ and (\ref{eq-Fgcl}) for $j=k+1$.

Finally, by definition we have $\partial _{k+1}(\widetilde{E_kK_k^{-1}})=
\partial _k(\tilde{F}_k)=1$. Hence
\begin{align*}
\partial _{k+1}(\widetilde{E_kK_k^{-1}}\gcl _{k+1}-q^{-2}\gcl _{k+1}
\widetilde{E_kK_k^{-1}})&=
\partial _k(\tilde{F}_k\gcl _k-q^2\gcl _k\tilde{F}_k)=2.
\end{align*}
This and Proposition \ref{s-clbasis} imply that (\ref{eq-Egcl}) and
(\ref{eq-Fgcl}) are fulfilled for $j=k+1$ and $j=k$, respectively.
\end{bew}

\begin{bew}[ of Theorem \ref{t-homuqgcl}]
We have to check that $\pi $ maps all relations
(\ref{eq-UK})--(\ref{eq-USerreF})
of $\Uqso N$ to zero. Instead of this we prove that the
corresponding equivalent relations, where the elements
$K_i,K_i^{-1},E_iK_i^{-1}$ and $F_i$ are involved, are mapped to
zero by $\pi $.

Since $\pi (K_i)$ and $\pi (K_j^{-1})$, $i,j=1,2,\ldots ,n$ are elements of
$\Clz $, by Proposition \ref{s-clzcomm} they commute.
Moreover, $z:=\pi (K_i)\pi (K_i^{-1})-1$ and $z':=\pi (K_i^{-1})\pi (K_i)-1$
are even elements of $\Cl $ commuting with all $\gcl _k$, $k=1,2,\ldots ,N$
by (\ref{eq-Kgcl}). {}From Proposition \ref{s-mucenter} we get
$z=\lambda z_{\mu ,0}=\lambda 1$ and $z'=\lambda 'z_{\mu ,0}=\lambda '1$,
where $\mu =(1,1,\ldots ,1)$, $\lambda ,\lambda '\in \comp $. 
Hence $(z-\lambda )\varphi ^1_N=(z'-\lambda ')\varphi ^1_N=0$.
On the other hand,
$\pi (K_i)\varphi ^1_N=q^{(2-\epsilon )\delta _{in}}\varphi ^1_N$ and
$\pi (K_i^{-1})\varphi ^1_N=q^{(\epsilon -2)\delta _{in}}\varphi ^1_N$
and therefore $\lambda \varphi ^1_N=\lambda '\varphi ^1_N=0$.
This means $\lambda =\lambda '=0$ and hence equation (\ref{eq-UK})
is compatible with $\pi $.

The compatibility of equations (\ref{eq-UKEKF}) (the first has to be
multiplied by $K_j^{-1}$) with
$\pi $ follows from (\ref{eq-Kgcl}) and (\ref{eq-alambda}).

Since $\partial _i(\pi (E_i))=\partial _i(\pi (E_iK_i^{-1}))=-1$, we have
$\partial _i(\pi (E_i)^2)=-2$ for all $i=1,2,\ldots ,n$. By Proposition
\ref{s-clbasis} this means that $\pi (E_i)^2=0$ for all $i=1,2,\ldots ,n$.
Further, $\partial _{i+1}(\pi (E_i))=1$ for $i\leq n-1$ and
$\partial _i(\pi (E_{i+1}))=\partial _{i+2}(\pi (E_i))=0$ for
$i\leq n-2$. Hence
$\partial _i(\pi (E_i)\pi (E_{i+1})\pi (E_i))=-2$ and
$\partial _{i+2}(\pi (E_{i+1})\pi (E_i)\pi (E_{i+1}))=2$ for $i\leq n-2$.
Therefore $\pi (E_i)\pi (E_{i+1})\pi (E_i)=
\pi (E_{i+1})\pi (E_i)\pi (E_{i+1})=0$ for $i\leq n-2$. Similarly,
\begin{align*}
\pi (E_{n-2})\pi (E_n)\pi (E_{n-2})=
\pi (E_n)\pi (E_{n-2})\pi (E_n)=0,\quad &N=2n,\\
\pi (E_{n-1})\pi (E_n)\pi (E_{n-1})=0,\quad &N=2n+1.
\end{align*}
This proves that (\ref{eq-USerreE}) is compatible with $\pi $ for $a_{ij}<0$.
By similar arguments we obtain the same for (\ref{eq-USerreF}).

Finally we have to prove compatibility of the relations (\ref{eq-UEF})
for all $i,j=1,2,\ldots ,n$
and (\ref{eq-USerreE}), (\ref{eq-USerreF}) for $a_{ij}=0$ with $\pi $.
In order to do this we recall from equations
(\ref{eq-Egcl})--(\ref{eq-Fngcl}) that there exist complex numbers
$\hat{a}_{kj}^l$, $\hat{b}_{kj}^l$, $k=1,2,\ldots ,n$, $j,l=1,2,\ldots ,N$,
such that
\begin{align*}
\pi (E_kK_k^{-1})\gcl _j=\lambda _{k,j}^{-1}\gcl _j\pi (E_kK_k^{-1})
+\hat{a}_{kj}^l\gcl _l,\qquad
\pi (F_k)\gcl _j=\lambda _{k,j}^{-1}\gcl _j\pi (F_k)+\hat{b}_{kj}^l\gcl _l
\end{align*}
for all $k=1,2,\ldots ,n$, $j=1,2,\ldots ,N$.
Then
\begin{align*}
&\left(q_i^{-a_{ik}}\pi (E_iK_i^{-1})\pi (F_k)-\pi (F_k)\pi (E_iK_i^{-1})-
\delta _{ik}\frac{1-\pi (K_i^{-1})^2}{q_i-q_i^{-1}}\right)\gcl _j=\\
&\quad =q_i^{-a_{ik}}\pi (E_iK_i^{-1})(\lambda _{k,j}^{-1}\gcl _j\pi (F_k)
+\hat{b}_{kj}^l\gcl _l)-\pi (F_k)(\lambda _{i,j}^{-1}\gcl _j\pi (E_iK_i^{-1})
+\hat{a}_{ij}^l\gcl _l)\\
&\quad \phantom{=}-\delta _{ik}\gcl _j\frac{1-\lambda _{i,j}^{-2}
\pi (K_i^{-1})^2}{q_i-q_i^{-1}}\\
&\quad =\lambda _{i,j}^{-1}\lambda _{k,j}^{-1}\gcl _j\left(
q_i^{-a_{ik}}\pi (E_iK_i^{-1})\pi (F_k)-\pi (F_k)\pi (E_iK_i^{-1})-
\delta _{ik}\frac{1-\pi (K_i^{-1})^2}{q_i-q_i^{-1}}\right)\\
&\quad \phantom{=}+(q_i^{-a_{ik}}\lambda _{i,l}^{-1}\hat{b}_{kj}^l-
\lambda _{i,j}^{-1}\hat{b}_{kj}^l)\gcl _l\pi (E_iK_i^{-1})
+(q_i^{-a_{ik}}\lambda _{k,j}^{-1}\hat{a}_{ij}^l-\lambda _{k,l}^{-1}
\hat{a}_{ij}^l)\gcl _l\pi (F_k)\\
&\quad \phantom{=}+\left(q_i^{-a_{ik}}\hat{b}_{kj}^l\hat{a}_{il}^m
-\hat{a}_{ij}^l\hat{b}_{kl}^m
+\delta _{ik}\delta _{jm}\frac{\lambda _{i,j}^{-1}\lambda _{k,j}^{-1}
-1}{q_i-q_i^{-1}}\right)\gcl _m
\end{align*}
for all $i,k=1,2,\ldots ,n$ and $j=1,2,\ldots ,N$. Observe that
$T_1(E_iK_i^{-1})\gcl _j=\hat{a}_{ij}^l\gcl _l$,
$T_1(F_k)\gcl _j=\hat{b}_{kj}^l\gcl _l$ and
$T_1(K_k)\gcl _j=\lambda _{k,j}\gcl _j$
for all $i,k=1,2,\ldots ,n$ and $j=1,2,\ldots ,N$.
Since
\begin{gather*}
T_1(q_i^{-a_{ik}}K_i^{-1}F_k-F_kK_i^{-1})\gcl _j=0,\quad
T_1(q^{-d_ia_{ik}}E_iK_i^{-1}K_k^{-1}-K_k^{-1}E_iK_i^{-1})\gcl _j=0,\\
T_1(q_i^{-a_{ik}}E_iK_i^{-1}F_k-F_kE_iK_i^{-1}+
\delta _{ik}(K_i^{-2}-1)/(q_i-q_i^{-1}))\gcl _j=0,
\end{gather*}
the last two lines in the above expression vanish.
Together with $\im \pi \subset \Cle $ we obtain
\begin{align*}
z_{ik}:=q_i^{-a_{ik}}\pi (E_iK_i^{-1})\pi (F_k)-\pi (F_k)\pi (E_iK_i^{-1})-
\delta _{ik}\frac{1-\pi (K_i^{-2})}{q_i-q_i^{-1}}\in \mc{Z}^{\mu ,0}_N,
\end{align*}
where $\mu =\mu _{(i,k)}=(\mu _1,\mu _2,\ldots ,\mu _n)$,
$\mu _l=\lambda _{i,l}^{-1}\lambda _{k,l}^{-1}$.
Similarly one can prove that
\begin{align*}
z'_{ik}&:=\pi (E_iK_i^{-1})\pi (E_kK_k^{-1})-\pi (E_kK_k^{-1})\pi (E_iK_i^{-1})
\in \mc{Z}^{\mu ,0}_N,\\
z''_{ik}&:=\pi (F_i)\pi (F_k)-\pi (F_k)\pi (F_i)\in \mc{Z}^{\mu ,0}_N
\end{align*}
for $a_{ik}=0$ and the same $\mu =\mu _{(i,k)}$.
Since $\tau (z_{ik})=z_{ki}$ and $\tau (z'_{ik})=z''_{ki}$, it suffices to
show that $z_{ik}=0$ for $i\leq k$ and $z'_{ik}=0$ for $i,k=1,2,\ldots ,n$,
$a_{ik}=0$.
By Proposition \ref{s-mucenter} we obtain $z_{ik},z'_{ik}
\in \comp z_{\mu ,0}$, $\mu =\mu _{(i,k)}$.
Obviously, $z_{\mu ,0}\varphi ^1_N=1$. Hence we only have to prove that
$z_{ik}\varphi ^1_N=0$ for $i\leq k$ and $z'_{ik}\varphi ^1_N=0$
for $i,k=1,2,\ldots ,n$, $a_{ik}=0$.
Since $\pi (E_iK_i^{-1})\varphi ^1_N=0$ for all $i=1,2,\ldots ,n$,
we get $z'_{ik}\varphi ^1_N=0$.
If $i\leq k<n$ then $\pi (K_i^{-1})\varphi ^1_N=\varphi ^1_N$ and
$\pi (F_k)=0$. Hence $z_{ik}\varphi ^1_N=0$.
If $i<k=n$ then $\partial _i(\pi (E_iK_i^{-1})\pi (F_n))
=-1+\delta ^i_{n-1}\delta _{\epsilon 0}$. Moreover, for $i=n-1,k=n$
we get $\partial _n(\pi (E_{n-1}K_{n-1}^{-1})\pi (F_n))=2$. Hence for $i<k=n$
all summands of $z_{ik}$ vanish after multiplication by $\varphi ^1_N$.
We are left with the case $i=k=n$, $N=2n+\epsilon $.
Then $\pi (E_nK_n^{-1})\varphi ^1_N=0$ and $\pi (K_n^{-1})\varphi ^1_N=
q^{\epsilon -2}\varphi ^1_N$.
For $N=2n$ we compute
\begin{align*}
z_{nn}\varphi ^1_{2n}&=\left(q^{-4}\pi (E_nK_n^{-1})
\frac{q^{-1}}{c^2\qp }\gcl _{n-1}\gcl _n-\frac{1-q^{-4}}{q^2-q^{-2}}\right)
\varphi ^1_{2n}\\
&=\left(\frac{q^{-5}}{c^2\qp }\frac{q^{-1}}{c^2\qp }\gcl _{n+1}\gcl _{n+2}
\gcl _{n-1}\gcl _n-q^{-2}\right)\varphi ^1_{2n}\\
&=\left(\frac{q^{-6}}{c^4\qp ^2}\gcl _{n+1}c^2q^{N-2n+3}\qp \gcl _n
-q^{-2}\right)\varphi ^1_{2n}\\
&=\left(\frac{q^{-3}}{c^2\qp }c^2q^{N-2n+1}\qp -q^{-2}\right)
\varphi ^1_{2n}=0.
\end{align*}
Similarly, for $N=2n+1$ we obtain
\begin{align*}
&z_{nn}\varphi ^1_{2n+1}=\left(q^{-2}\pi (E_nK_n^{-1})
\frac{q^{-1}q'}{c^2\qp }\gcl _n\gcl _{n+1}-\frac{1-q^{-2}}{q-q^{-1}}\right)
\varphi ^1_{2n+1}\\
&\quad =\left(\frac{q^{-3}}{cq'}\pi (E_nK_n^{-1})\gcl _n-q^{-1}\right)
\varphi ^1_{2n+1}
=\left(\frac{q^{-3}}{cq'}\frac{q'}{c^2\qp }\gcl _{n+1}\gcl _{n+2}
\gcl _n-q^{-1}\right)\varphi ^1_{2n+1}\\
&\quad =\left(\frac{q^{-3}}{c^3\qp }\gcl _{n+1}c^2q^{N-2n+1}\qp -q^{-1}\right)
\varphi ^1_{2n+1}
=\left(\frac{q^{-1}}{c}\gcl _{n+1}-q^{-1}\right)\varphi ^1_{2n+1}=0.
\end{align*}
\end{bew}

\begin{bew}[ of Theorem \ref{t-uqgreps}]
First let us prove (i). Since $\pi _\mathrm{ad}$ is a representation of
$\Uqso N$ on $\Cl $, it is sufficient to check the invariance of $V$
under the action of the generators $f=K_i,K_i^{-1},E_iK_i^{-1},F_i$
of $\Uqso N$. By the definition of $\Uqso N$, formula (\ref{eq-piad})
for these $f$ takes the form
\begin{gather*}
(\pi _\mathrm{ad}K_i)(\gcl _j)=\pi (K_i)\gcl _j\pi(K_i^{-1}),\quad
(\pi _\mathrm{ad}K_i^{-1})(\gcl _j)=\pi (K_i^{-1})\gcl _j\pi (K_i),\\
(\pi _\mathrm{ad}E_iK_i^{-1})(\gcl _j)=\pi (E_iK_i^{-1})\gcl _j
-\pi (K_i^{-1})\gcl_j\pi (K_i)\pi (E_iK_i^{-1}),\\
(\pi _\mathrm{ad}F_i)(\gcl _j)=\pi (F_i)\gcl _j-\pi (K_i^{-1})\gcl_j\pi (K_i)
\pi (F_i)
\end{gather*}
for $i=1,2,\ldots ,n$ and $j=1,2,\ldots ,N$.
{}From this and Lemma \ref{l-pirels} the invariance of $V$ under the
action $\pi _\mathrm{ad}$ of $\Uqso N$ immediately follows.
Moreover,
\begin{align*}
(\pi _\mathrm{ad}K_i)(\gcl _N)&=\lambda _{i,N}\gcl _N=q^{2\delta _{i1}}\gcl _N,
\\
(\pi _\mathrm{ad}E_i)(\gcl _N)&=
(\pi _\mathrm{ad}E_iK_i^{-1})((\pi _\mathrm{ad}K_i)(\gcl _N))=
q^{2\delta _{i1}}(\pi _\mathrm{ad}E_iK_i^{-1})(\gcl _N)=0,
\end{align*}
and therefore the $N$-dimensional representation $\pi _\mathrm{ad}$
of $\Uqso N$ is isomorphic to the vector representation
of $\Uqso N$ with higest weight $(1+\delta _{N3},0,0,\ldots ,0)$ and
highest weight vector $\gcl _N$.

Now we turn to (ii). Since $\mc{I}^\nu _\mathrm{L}$, $\nu \in \{+,-\}$,
is a left ideal of $\Cl[2n{+}1]$, the mapping $\pi ^\nu $ gives a representation
of $\Uqso{2n+1}$. The dimension of this representation is
$\dim \mc{I}^\nu _\mathrm{L}=2^n$ by Proposition \ref{s-idealbasis}.
On the other hand, from $\pi (E_iK_i^{-1})\varphi ^\nu _{2n+1}=0$ and
$\pi (K_i)\varphi ^\nu _{2n+1}=q^{\delta _{in}}$ for all $i=1,2,\ldots ,n$
we obtain that the $2^n$-dimensional representation of $\Uqso{2n+1}$
with highest weight $(0,0,\ldots ,0,1)$ and highest weight vector
$\varphi ^\nu _{2n+1}$ is a subrepresentation of $\pi ^\nu $.
{}From this the assertion follows.

Finally we prove (iii). Since $\pi (f)\in \Cle[2n]$
for all $f=K_i,K_i^{-1},E_iK_i^{-1},F_i$, we have $\pi (\Uqso{2n})
\subset \Cle[2n]$. Moreover, $\varphi ^1_{2n}$ is homogeneous with respect
to the $\mathbb{Z}_2$-grading $\partial _0$ of $\Cl[2n]$, and hence
$\Cls[2n]{\nu }\varphi ^1_{2n}$, $\nu =+,-$, is an invariant subspace
under left multiplication by $\Cle[2n]$, and therefore by $\Uqso{2n}$ as well.
By Proposition \ref{s-idealbasis}
the dimension of the vector space $\Cls[2n]{\nu }\varphi ^1_{2n}$, $\nu =+,-$,
is $2^{n-1}$. Again we obtain
\begin{align*}
(\pi ^+K_i)(\varphi ^1_{2n})&=q^{2\delta _{in}}\varphi ^1_{2n},&
(\pi ^+E_i)(\varphi ^1_{2n})&=\pi (E_i)\varphi ^1_{2n}=
\pi (E_iK_i^{-1})q^{2\delta _{in}}\varphi ^1_{2n}=0
\end{align*}
for $i=1,2,\ldots ,n$.
Hence $\pi ^+$ contains as a subrepresentation the $2^{n-1}$-dimensional
highest weight representation with highest weight $(0,\ldots ,0,1)$ and
highest weight vector $\varphi ^1_{2n}$.
Further,
\begin{align*}
(\pi ^-K_i)(\gcl _n\varphi ^1_{2n})&=q^{2\delta _{in}}
\left(1+\frac{\lambda _{i,n}-1}{c^2\qp q^{N+1-2n}}\gcl _n\gcl _{n+1}\right)
\gcl _n\varphi ^1_{2n}\\
&=q^{2\delta _{in}}\gcl _n\left(1+\frac{\lambda _{i,n}-1}{c^2\qp q^{N+1-2n}}
c^2q^{N-2n+1}\qp \right)\varphi ^1_{2n}\\
&=q^{2\delta _{in}}\lambda _{i,n}\gcl _n\varphi ^1_{2n}
=q^{2\delta ^i_{n-1}}\gcl _n\varphi ^1_{2n},\\
(\pi ^-E_iK_i^{-1})(\gcl _n\varphi ^1_{2n})&=\pi (E_iK_i^{-1})\gcl _n
\varphi ^1_{2n}
=\delta _{in}\frac{q^{-1}}{c^2\qp }\gcl _{n+1}\gcl _{n+2}\gcl _n
\varphi ^1_{2n}=0
\end{align*}
for $i=1,2,\ldots ,n$.
Hence $\pi ^-$ contains as a subrepresentation the $2^{n-1}$-dimensional
highest weight representation with highest weight $(0,0,\ldots ,0,1,0)$ and
highest weight vector $\gcl _n\varphi ^1_{2n}$. This proves the theorem.
\end{bew}

\bibliographystyle{mybib}
\bibliography{quantum}

\end{document}